\documentclass[psamsfonts,11pt]{amsart}
\usepackage{amsmath,amssymb,amsthm,amscd} 
\usepackage[left=3.5cm,right=3.5cm,top=3.5cm,bottom=3.5cm]{geometry}
\usepackage{amssymb,fge} 
\usepackage{tikz-cd}
\usepackage{mathtools}
\usepackage{hyperref}  
\usepackage{comment}
\usepackage{graphicx}
\usepackage{dirtytalk}

\newtheorem{thm}{Theorem}
\newtheorem{cor}[thm]{Corollary}
\newtheorem{prop}[thm]{Proposition}
\newtheorem{lem}[thm]{Lemma}

\newtheorem{defn}[thm]{Definition}
\newtheorem{question}[thm]{Question}
\newtheorem{conjecture}[thm]{Conjecture}

\newcommand\scalemath[2]{\scalebox{#1}{\mbox{\ensuremath{\displaystyle #2}}}}

\theoremstyle{definition} 

\newtheorem{rem}[thm]{Remark}

\theoremstyle{remark}

\numberwithin{thm}{section}
\numberwithin{equation}{thm}
\newcommand{\be}{\begin{equation}}
\newcommand{\ee}{\end{equation}}

\DeclareMathOperator{\Prep}{Preper}

\DeclareMathOperator{\rank}{rank}

\title[]{On the proportion of irreducible polynomials in unicritically generated semigroups}
\author[]{Wade Hindes}
\author[]{Reiyah Jacobs}
\author[]{Benjamin Keller}
\author[]{Albert Kim}
\author[]{Peter Ye}
\author[]{Aaron Zhou}

\begin{document}
\maketitle
\begin{abstract} Let $p$ be a prime number and let $S=\{x^p+c_1,\dots,x^p+c_r\}$ be a finite set of unicritical polynomials for some $c_1,\dots,c_r\in\mathbb{Z}$. Moreover, assume that $S$ contains at least one irreducible polynomial over $\mathbb{Q}$. Then we construct a large, explicit subset of irreducible polynomials within the semigroup generated by $S$ under composition; in fact, we show that this subset has positive asymptotic density within the full semigroup when we count polynomials by degree. In addition, when $p=2$ or $3$ we construct an infinite family of semigroups that break the local-global principle for irreducibility. To do this, we use a mix of algebraic and arithmetic techniques and results, including Runge's method, the elliptic curve Chabauty method, and Fermat's Last Theorem.  
\end{abstract}
\section{Introduction}
Given a field $K$, it is often an interesting and difficult problem to explicitly construct irreducible polynomials over $K$ of arbitrarily large degree. For instance when $K/\mathbb{Q}$ is a number field, there are relatively few known tests for irreducibilty (e.g., reduction modulo primes, Eisenstein's criterion, etc.), all of which exploit some particular arithmetic aspect of the given polynomial. This is true despite the fact that \say{most} polynomials of a given degree are irreducible by Hilbert's Irreducibility Theorem. In this paper, we study the problem of constructing irreducible polynomials of arbitrarily large degree via composition. In particular, we hope to construct a \say{large} set of irreducible polynomials in the following sense. 
\begin{defn}\label{def:proportion} Let $K$ be a field, let $S=\{f_1,\dots,f_r\}$ be a finite set of elements of $K[x]$, and let $M_S$ be the semigroup generated by $S$ under composition. Then we say that $M_S$ contains a positive proportion of irreducible polynomials over $K$ if:
\[\liminf_{B\rightarrow\infty}\,\frac{\#\big\{f\in M_S\,: \deg(f)\leq B\;\text{and $f$ is irreducible in $K[x]$}\big\}}{\#\big\{f\in M_S\,: \deg(f)\leq B\big\}}>0.\]
\end{defn}
With a little thought one can see that the generating set $S$ should contain at least one irreducible polynomial, otherwise the proportion above is zero. But is this enough? What if we assume also that the semigroup $M_S$ is sufficiently large, for instance, free? As motivation for future work on this problem, we pose the following question.
\begin{question}\label{quest} Let $K/\mathbb{Q}$ be a number field and let $S=\{f_1,\dots,f_r\}$ be a finite set of elements of $K[x]$ with $r\geq2$. If $S$ contains at least one irreducible polynomial and $M_S$ is a free semigroup, does $M_S$ contain a positive proportion of irreducible polynomials? 
\end{question}
Perhaps, the first non-trivial case to investigate this question is to consider sets $S$ of quadratic polynomials with integer coefficients. Indeed, using a mix of algebraic and arithmetic techniques, we answer the question above in the affirmative in this case. To state our results, we make the following definitions. 
\begin{defn}\label{def:types} \textup{Let $p$ be a prime. We say that $\phi(x)=x^p+c$ is of Type (I) if $c=s^p-s^{p^2}$ for some $s\in\mathbb{Q}$. Equivalently, $\phi$ is of Type (I) if it has a rational fixed point that is $p$th power. Likewise, for $p=2$ we say that $\phi=x^2+c$ is of Type (II) if $c=-1-s^2-s^4$ for some $s\in\mathbb{Q}$. Equivalently, $\phi=x^2+c$ is of Type (II) if it has a rational point of exact period two that is a square}.
\end{defn}
With the notions above, our first main result is the following. Here and throughout, irreducibility refers to irreducibility over the field of rational numbers. 
\begin{thm}\label{thm:main} Let $S=\{x^2+c_1,\dots,x^2+c_r\}$ for some distinct $c_1,\dots,c_r\in\mathbb{Z}$ and assume that $S$ contains at least one irreducible polynomial. Then the following statements hold. 
\begin{enumerate}
\item If $r=1$, then $M_S$ is a set of irreducible polynomials. \vspace{.2cm}
\item If $S$ contains at least one irreducible polynomial $\phi$ that is neither of Type (I) nor of Type(II), then 
    \[\{\phi^4\circ f\,:\, f\in M_S\}\] is a set of irreducible polynomials. \vspace{.2cm}
\item If $S$ contains two irreducible polynomials $\phi_1$ and $\phi_2$ that are both of Type (I) or of Type(II), though not necessarily of the same type, then either
\begin{equation*}
\{\phi_1^3\circ\phi_2\circ f\, :\, f\in M_S\}
\quad
\text{or}
\quad
\{\phi_2^3\circ\phi_1\circ f\, :\, f\in M_S\}
\end{equation*}
is a set of irreducible polynomials. \vspace{.2cm}
\item If $S$ contains an irreducible $\phi_1(x)=x^2+s^2-s^4$ of Type (I) and a reducible $\phi_2(x)=x^2-t^2$ such that $s^2\neq t^2$, then 
\[
\{\phi_1^3\circ\phi_2\circ\phi_1\circ f\, :\, f\in M_S\}
\]
is a set of irreducible polynomials.\vspace{.2cm}
\item If $S$ contains an irreducible $\phi_1(x)=x^2+s^2-s^4$ of Type (I) and a reducible $\phi_2(x)=x^2-s^2$, then 
\[
\{\phi_1^2\circ\phi_2\circ\phi_1\circ f\, :\, f\in M_S\}
\]
is a set of irreducible polynomials.\vspace{.2cm}
\item If $S$ contains an irreducible polynomial $\phi_1$ of Type(II) and a reducible polynomial $\phi_2$, then  
\[
\{\phi_1^3\circ\phi_2\circ\phi_1\circ f\, :\, f\in M_S\}
\]
is a set of irreducible polynomials. \vspace{.2cm}
\end{enumerate}
In particular, \vspace{.15cm}
\[\liminf_{B\rightarrow\infty}\,\frac{\#\big\{f\in M_S\,: \deg(f)\leq B\;\text{and $f$ is irreducible}\big\}}{\#\big\{f\in M_S\,: \deg(f)\leq B\big\}}\geq \frac{1}{(\#S)^5}. \vspace{.15cm}\]
Hence, $M_S$ contains a positive proportion of irreducible polynomials if and only if it contains at least one irreducible polynomial. 
\end{thm}
In particular, motivated by the case of quadratic polynomials with integer coefficients above, we make the following conjecture. 

\begin{conjecture} Let $S=\{x^2+c_1,\dots,x^2+c_r\}$ for some distinct $c_1,\dots,c_r\in\mathbb{Q}$ and $r\geq2$. Then $M_S$ contains a positive proportion of irreducible polynomials if and only if it contains at least one irreducible polynomial. 
\end{conjecture}
\begin{rem} Assuming Poonen's uniformity conjecture \cite{Poonen}, we prove the conjecture above for sets $S$ containing at least one \emph{stable} polynomial \cite[\S5]{jones2013galois} that is neither of Type (I) or Type (II) by adapting our methods in Section 3 and encorporating height bounds for rational points on curves; see Theorem  \ref{thm:rational} below.      
\end{rem}
Likewise, we prove a similar result for sets of the form $S=\{x^p+c_1,\dots,x^p+c_r\}$ when $p$ is an odd prime and $c_1,\dots,c_r\in\mathbb{Z}$. However, there is a catch: when $p>3$ we must assume that $S$ \emph{does not} take the very special form, $S=\{x^p+t^p-t^{p^2},x^p+t^p\}$ for some $t\in\mathbb{Z}$. Indeed, this case remains open. On the other hand, outside of this case we use a mix of algebraic and arithmetic techniques and results, including Runge's method \cite[\S5]{schoof2010catalan}, the elliptic curve Chabauty method \cite{Bruin}, and Fermat's Last Theorem \cite{FLT}, to prove the following result; see Theorem \ref{thm:oddprime,explicit} below for an explicit version.     

\begin{thm}\label{thm:intro-oddprime}
Let $S=\{x^p+c_1,\dots,x^p+c_r\}$ for some odd prime $p$ and some distinct $c_1,\dots,c_r\in\mathbb{Z}$. Moreover, assume that $S$ contains at least one irreducible polynomial. Then the following statements hold: \vspace{.1cm}   
\begin{enumerate}
    \item If $p=3$, then $M_S$ contains a positive proportion of irreducible polynomials. \vspace{.1cm} 
    \item If $p>3$ and $S\neq\{x^p+t^p-t^{p^2},x^p+t^p\}$ for all $t\in\mathbb{Z}$, then $M_S$ contains a positive proportion of irreducible polynomials.  
\end{enumerate}
\end{thm}

The special polynomials in Theorem \ref{thm:intro-oddprime} are interesting for an additional reason: they are reducible modulo every prime. Moreover, the construction of irreducible polynomials with this property have a long history (going back at least to David Hilbert, who pointed out that $x^4+1$ has this property). With this in mind, we make the following definition. 
\begin{defn} A monic $f\in\mathbb{Z}[x]$ breaks the local-global principle for irreducibility if $f$ is irreducible in $\mathbb{Q}[x]$ and its image $\overline{f}\in\mathbb{F}_q[x]$ is reducible for every prime $q$. 
\end{defn}
In particular, we show that the special polynomials in Theorem \ref{thm:intro-oddprime} and in case (4) of Theorem \ref{thm:main} break the local-global principle for irreduciblity; compare to \cite{jones2012iterative}.  
\begin{thm} \label{thm:local-global} Let $p=2$ or $p=3$, let $t\in\mathbb{Z}\setminus\{0,\pm{1}\}$, and let  
\[S=\big\{x^p+t^p-t^{p^2},\;x^p+(-1)^{p-1}t^p,\;x^p+c_1,\;\dots,\;x^p+c_r\big\}\] 
for some arbitrary $c_1,\dots,c_r\in\mathbb{Z}$ (including the case of no extra $c$'s). Then:\vspace{.15cm}  
\begin{enumerate}
\item If $p=2$, then $\{\phi_1^2\circ\phi_2\circ\phi_1\circ f\, :\, f\in M_{S}\}$
is a set of polynomials that break the local-global principal for irreducibility. \vspace{.2cm} 
\item If $p=3$, then $\{\phi_1\circ\phi_2\circ\phi_1\circ f\, :\, f\in M_{S}\}$
is a set of polynomials that break the local-global principal for irreducibility. \vspace{.15cm}  
\end{enumerate}
In particular, the proportion of polynomials in $M_S$ that break the local-global principal for irreducibility is at least ${(\#S)}^{-4}>0$.  
\end{thm}
Finally, we adapt ideas used to prove the theorems above to establish the following partial result for unicritical sets with \emph{rational} coefficients. In what follows, $\phi\in\mathbb{Q}[x]$ is called \emph{stable} if $\phi^n$ is irreducible in $\mathbb{Q}[x]$ for all $n\geq1$; see \cite[\S5]{jones2013galois}.    
\begin{thm}\label{thm:rational} Let $S=\{x^p+c_1,\dots,x^p+c_r\}$ for some prime $p$ and some distinct $c_1,\dots,c_r\in\mathbb{Q}$. Moreover, assume the following:\vspace{.1cm}  
\begin{enumerate}
    \item If $p=2$, assume that $S$ contains a stable polynomial $\phi$ that is neither of Type (I) nor of Type (II) and that Poonen's Conjecture \cite{Poonen} holds.\vspace{.2cm}
    \item If $p>2$, assume that $S$ contains a stable polynomial $\phi$ that is not of Type (I). \vspace{.2cm}   
\end{enumerate}
Then there exists $n=n_\phi\gg0$ such that 
\[\{\phi^{n}\circ f\,:\, f\in M_S\}\]
is a set of irreducible polynomials. In particular, $M_S$ contains a positive proportion of irreducible polynomials in this case.  
\end{thm}
An outline of our paper is as follows. In Section 2, we first review an irreducibility test for compositions of unicritical polynomials that relates powers in dynamical orbits to reducibility. We then prove two crucial classification theorems, both of which say that if $\phi=x^p+c$ is an irreducible, unicritical polynomial defined over $\mathbb{Z}$ and $\phi$ possesses a large iterate with a $p$th powered image, then $\phi$ is necessarily of Type (I) or of Type (II). From here, we construct irreducible polynomials in the semigroups generated by sets of the form $S=\{x^p+c_1,\dots,x^p+c_r\}$ for $p=2$ in Section 3 and for odd $p$ in Section 4. Moreover, we include the local-global result in Section 4 as well. Lastly, we prove Theorem \ref{thm:rational} in Section 5. For some {\tt{Magma}} \cite{Magma} code accompanying this work, see  \href{https://sites.google.com/a/alumni.brown.edu/whindes/research/code}{Mathworks 2023}.    
\\[7pt] 
\noindent\textbf{Acknowledgements:} We thank the Mathworks program at Texas State University for supporting this research. We also thank Rafe Jones and Mohammad Sadek for some useful conversations related to the work in this paper. 
\section{Powers in Orbits and Irreducibility}
We begin with the following irreducibility test for iterates in unicritically generated semigroups; see \cite[Theorem 8]{JonesIMRN} for more general statement and proof. 
\begin{prop}\label{prop:stability} Let $p$ be a prime, let $K$ be a field of characteristic not equal to $p$, let $f(x)=x^p+c$ for some $c\in K$, and let $g(x)\in K[x]$ be a monic, irreducible, and separable polynomial. Moreover, when $p=2$ we assume that the degree of $g$ is even. Then, if $g(f(0))$ is not a $p$th power in $K$, we have that $g(f(x))$ is irreducible and separable in $K[x]$. Furthermore, this statement becomes an \say{if and only if} statement when $K$ has the following property: for every finite extension $L/K$, the induced norm homomorphism $N_{L/K}: L^*/(L^*)^p\rightarrow K^*/(K^*)^p$ is injective. 
\end{prop}
\begin{rem}\label{rem:finitefield} In particular, the result above is an \say{if and only if} statement when $K$ is a finite field of characteristic coprime to $p$.     
\end{rem}
Now that we have a link between powers and reducibility in unicritically generated semigroups, we next classify the integral quadratic polynomials possessing a large iterate with a square image. In particular, the only such polynomials must have a square fixed-point or a $2$-cycle containing a square (i.e., be of Type I or Type II).
\begin{thm}\label{thm:4th=square} 
Let $\phi(x)=x^2+c$ for some $c\in\mathbb{Z}$ and assume that $\phi$ is irreducible over $\mathbb{Q}$. Then, if $y^2=\phi^4(a)$ for some $a,y\in\mathbb{Z}$, there exists $s\in\mathbb{Z}$ such that one of the following statements must hold. \vspace{.15cm} 
\begin{enumerate}
    \item $c=s^2-s^4$ and $a=\pm{s^2}$. In particular, $\phi$ is of Type (I). \vspace{.25cm} 
     \item $c=-1-s^2-s^4$ and $a=\pm{s^2}$. In particular, $\phi$ is of Type (II).
\end{enumerate}
\end{thm} 
\begin{proof} Let $\phi=x^2+c$ for some non-zero $c\in\mathbb{Z}$. Then we begin with a simple observation (that enables us to sidestep the use of elliptic curves as in \cite{Mathworks}):
\\[5pt]
\textbf{Fact 1:} If $y^2=\phi(a)$ for some $a,y\in\mathbb{Z}$, then $|a|\leq |c|$.\\[5pt]
Indeed, if $\phi(a)=y^2$ then $(a-y)(a+y)=-c$ and $a\pm y$ are non-zero integer factors of $c$. Hence, $|a\pm y|\leq |c|$ and $|a|\leq |c|$ follows from the triangle inequality.   
In particular, we note that for positive $c>0$ and all $m\geq2$: 
\begin{equation}\label{eq:c>0}
\phi^m(x)\geq\phi^2(x)=x^4+2x^2c+c^2+c\geq c^2+c> c.
\end{equation} 
Hence, $\phi^n(x)\neq y^2$ for all $x,y\in\mathbb{Z}$ and $n\geq3$, which may be seen by substituting $a=\phi^{n-1}(x)$ and $m=n-1$ into Fact (1) and \eqref{eq:c>0}. In particular, we may assume that $c<-1$ (since $x^2-1$ is reducible). In this case, we note that since $-c$ is not a square there exists a unique, positive $s\in\mathbb{Z}$ such that $s^2<-c<(s+1)^2$. In particular, $-c=s^2+r$ some remainder term $0<r<2s+1$, and for the rest of the proof we write
\begin{equation}\label{boxed}
\boxed{\phi=x^2-s^2-r\;\;\;\text{for some $s\geq1$ and $1\leq r<2s+1$.}}
\end{equation} 
From here, we make a few more observations.\\[5pt] 
\textbf{Fact 2:} If $a\in\mathbb{Z}$ is such that $|\phi(a)|< 2s$, then $a=\pm{s},\pm{(s+1)}$. \\[5pt]
Indeed, if $|a|\leq s-1$, then 
\[\phi(a)=\phi(|a|)\leq \phi(s-1)=(s-1)^2-s^2-r\leq-2s+1-r\leq -2s.\]
Hence, $|\phi(a)|\geq 2s$ in this case. Likewise, if $|a|\geq s+2$, then 
\[\phi(a)=\phi(|a|)\geq \phi(s+2)=(s+2)^2-s^2-r\geq 2s+3>2s.\]
Therefore, $|\phi(a)|>2s$ in this case. In particular, $|\phi(a)|< 2s$ implies that $|a|=s,s+1$ as claimed. 
Next, we note the following fact.  
\\[5pt] 
\textbf{Fact 3:} If $|a|\geq 2s$ and $s\geq3$, then $|\phi^n(a)|>|c|$ for all $n\geq1$.\\[5pt] 
Indeed, if $|a|\geq 2s$ and $s\geq3$, then \[\phi(a)=\phi(|a|)\geq\phi(2s)\geq 3s^2-2s-1> s^2+2s+1>s^2+r=|c|.\] 
And $|\phi^n(a)|=\phi^n(a)>|c|$ follows easily by induction and the fact that $\phi(|c|)=c^2+c\geq|c|$ for all $|c|>1$.
Finally, we have the following fact. \\[5pt]  
\textbf{Fact 4:} If $s\geq3$ and $\phi^n(a)=y^2$ for some $n\geq3$ and $a,y\in\mathbb{Z}$, then $a=\pm{s}$ or $a=\pm{(s+1)}$.\\[5pt] 
To see this, note that if $a\not\in\{\pm{s},\pm{(s+1)}\}$, then $|\phi(a)|\geq2s$ by Fact $1$. In particular, replacing $a$ with $\phi(a)$ in Fact $3$ we see that $|\phi^m(a)|>c$ for all $m\geq2$ for such $a$. But if $y^2=\phi^n(a)=\phi(\phi^{n-1}(a))$, then Fact $1$ implies that $|\phi^{n-1}(a)|\leq c$. Therefore, it must be the case that $n\leq2$ whenever $a\not\in\{\pm{s},\pm{(s+1)}\}$ as claimed. 
\\[6pt]      
\indent Finally, we have the tools in place to complete the proof of the Theorem \ref{thm:4th=square}. Assume first that $s\geq3$ and that $\phi^4(a)=y^2$ for some integers $a$ and $y$. Without loss, we may assume that $a\geq0$. Then Fact 4 implies $a=s$ or $a=s+1$. From here, we proceed in cases. \\[3pt]
Case (1): Suppose first that $a=s$. Then $y^2=\phi^4(a)=\phi^3(\phi(a))=\phi^3(-r)$ and Fact 3 implies that $-r=-s,-(s+1)$. However, when $-r=-s$ we see that $-s=\phi(s)=\phi(-s)$. In particular, $y^2=\phi^4(s)=\phi^3(-s)=-s$, which contradicts the assumption that $s$ is positive. Therefore, $-r=-(s+1)$ and 
\[\phi^2(s)=\phi(\phi(s))=\phi(-(s+1))=(-(s+1))^2-s^2-(s+1)=s.\]
But then $y^2=\phi^4(s)=s$ and $-r=-(s+1)=-(y^2+1)$. Hence,  
\[\text{$\phi=x^2-s^2-r=x^2-y^4-y^2-1$\,\; and \;\,$a=s=y^2$,}\] 
and we obtain a description of $\phi$ and $a$ of the form in statement (2) of Theorem \ref{thm:4th=square} by replacing $s$ with $y$.
\\[3pt]
Case (2): Suppose next that $a=s+1$. Then since $y^2=\phi^3(\phi(a))$, we see that $\phi(a)\in\{\pm{s},\pm{(s+1)}\}$ by Fact 4. On the other hand, \[\phi(a)=\phi(s+1)=(s+1)^2-s^2-r=2s+1-r\geq0\] 
by \eqref{boxed}. Hence, $\phi(a)=2s+1-r=s$ or $\phi(a)=2s+1-r=s+1$ Suppose first that $\phi(a)=2s+1-r=s$. Then 
\[\phi^2(-(s+1))=\phi^2(s+1)=\phi^2(a)=\phi(s)=s^2-s^2-r=-r=-(s+1).\]
In particular, it follows that $y^2=\phi^4(a)=-(s+1)$, a contradiction. Therefore, we may assume that $\phi(a)=2s+1-r=s+1$. But then $\phi(s+1)=s+1$ and $y^2=\phi^4(a)=\phi^4(s+1)=s+1.$
In particular, 
\[\phi=x^2-s^2-r=x^2-(1-y^2)^2-(y^2-1)=x^2+y^2-y^4\]
and $a=y^2$, and we obtain a description of $\phi$ and $a$ of the form in statement (1) of Theorem \ref{thm:4th=square} by replacing $s$ with $y$.

Therefore, it suffices to verify the theorem for irreducible $\phi_c=x^2+c$ of the form in \eqref{boxed} for $s=1,2$. In particular, it remains to consider $c=\{-2,-3,-5,-6,-7,-8\}$. However, only $c=-3$ yields a solutions $(a,y)$ to the equation $\phi_c^4(a)=y^2$ modulo $24$. Moreover, the theorem is easily verified in this case: $c=-3=-1-1^2-1^4$ and $a=\pm{1}$, fitting form (2) in the statement of the theorem.  
%> for c in [-2,-3,-5,-6,-7,-8] do 
%> for a,y in Integers(24) do 
%> if (((a^2+c)^2+c)^2+c)^2+c eq y^2 then 
%> print c; 
%> end if; 
%> end for; 
%> end for; 
\end{proof}

In particular, we obtain the following result for sets of quadratic polynomials of the form $x^2+c$ containing at least one irreducible polynomial that is neither of Type (I) nor of Type (II). 

\begin{cor}\label{cor:irre&non-special} Let $S=\{x^2+c_1,\dots,x^2+c_r\}$ for some $c_i\in\mathbb{Z}$. Moreover, assume that $S$ contains at least one irreducible polynomial $\phi$ that is neither of Type (I) nor of Type (II). Then $\{\phi^4\circ f\,:\, f\in M_S\}$ is a set of irreducible polynomials in $\mathbb{Q}[x].$
\end{cor}
\begin{proof} First we note that \cite[Corollary 1.3]{Stoll} implies that $ \phi^4(x)$ is irreducible. Therefore, \cite[Proposition 2.1]{Mathworks} implies that if a polynomial of the form $ \phi^4\circ f$ for some $f\in M_S$ is reducible, then $\phi^4(a)=y^2$ for some $a,y\in\mathbb{Z}$. However, in this case it follows from Theorem \ref{thm:4th=square} that $\phi$ is of Type (I) or of Type (II). The claim follows.   
\end{proof}

On the other hand, building on the proof of Theorem \ref{thm:4th=square}, we have the following refinement about squares in orbits of polynomials of Type (I). 
\begin{lem}\label{lem:family1refinement} Let $\phi(x)=x^2+s^2-s^4$ be of Type (I) for some $s\in\mathbb{Z}\setminus\{0,\pm{1}\}$ and  assume that $\phi^3(a)=y^2$ for some $a,y\in\mathbb{Z}$. Then $a=\pm{s^2}$. 
\end{lem}
\begin{proof} Assume that $\phi^3(a)=y^2$ and $a\geq0$. Then $(\phi^2(a)-y)\cdot(\phi^2(a)-y)=s^2-s^4$ and $|\phi^2(a)|\le s^2+s^4$ follows from the triangle inequality. However, we note that if $a\geq s^2+1$, then    \[\phi^2(a)\geq\phi^2(s^2+1)=8s^4 + 7s^2 + 1>s^2+s^4.\]  
Hence, $\phi^3(a)\neq y^2$ in this case. Likewise, if $a\leq s^2-2$, then 
\[\phi(a)\leq \phi(s^2-1)=-3s^2+4<0.\]
Hence, $|\phi(a)|=-\phi(a)\geq 3s^2-4$ and therefore 
\[\phi^2(a)=\phi(|\phi(a)|)\geq \phi(3s^2-4)=8s^4 - 23s^2 + 16>s^2+s^4;\]
here we use that $|s|\geq2$. In particular, $\phi^3(a)\neq y^2$ in this case also, and we deduce that $a=s^2$ or $a=s^2-1$. On the other hand, when $a=s^2-1$ we have that $\phi^3(a)=1-s^2<0$ and so $\phi^3(a)$ cannot be a square. Hence, $a=s^2$ as desired. 
\end{proof}

Likewise, we have the following refinement about squares in orbits of polynomials of Type (II). 
\begin{lem}\label{lem:family2refinement} Let $\phi(x)=x^2-1-s^2-s^4$ be of Type (II) for some non-zero $s\in\mathbb{Z}$ and assume that $\phi^3(a)=y^2$ for some $a,y\in\mathbb{Z}$. Then $a=\pm{(s^2+1)}$. 
\end{lem}
\begin{proof} Assume that $\phi^3(a)=y^2$ and $a\geq0$. Then $(\phi^2(a)-y)\cdot(\phi^2(a)-y)=1+s^2+s^4$ and $|\phi^2(a)|\leq 1+s^2+s^4$ follows from the triangle inequality. However, we note that if $a\geq s^2+2$, then    \[\phi^2(a)\geq\phi^2(s+2)=8s^4 + 17s^2 + 8>1+s^2+s^4.\] 
Hence, $\phi^3(a)\neq y^2$ in this case. Likewise, if $a\leq s^2-1$, then 
\[\phi(a)\leq \phi(s^2-1)=-3s^2.\]
Hence, $|\phi(a)|\geq 3s^2$ and therefore 
\[\phi^2(a)=\phi(|\phi(a)|)\geq \phi(3s^2)=8s^4 - s^2 - 1>1+s^2+s^4;\]
here we use that $|s|\geq1$. In particular, $\phi^3(a)\neq y^2$ in this case also, and we deduce that $a=s^2$ or $a=s^2+1$. However, $\phi^3(s^2)=-(1+s^2)<0$ (so cannot be a square), and thus $a=s^2+1$ as claimed.     
\end{proof}

Finally, we end this section by proving an analog of Theorem \ref{thm:4th=square} for odd primes $p$.   
\begin{thm}\label{thm:pthpower} Let $\phi(x)=x^p+c$ for some odd prime $p$ and $c\in\mathbb{Z}$. Moreover, assume that $\phi$ is irreducible. If $\phi^3(a)=y^p$ for some $a,y\in\mathbb{Z}$, then $c=s^p-s^{p^2}$ and $a=s^p$ for some $s\in\mathbb{Z}$. Equivalently, $\phi$ is of Type (I) and $a$ is a pth powered fixed point of $\phi$.  
\end{thm}
We first establish a bound that will greatly help us in our proof; compare to Fact 1 in the proof of Theorem \ref{thm:4th=square}.  
\begin{lem}\label{lem:oddpbound} 
If $x^p + c = y^p$ for some $x,y,c \in \mathbb{Z}$ and $c \neq 0$, then $|x| \leq \left(\frac{|c|}{p}\right)^{\frac{1}{p-1}} + 1$. 
\end{lem}
\begin{proof}  
To see this, first note that by replacing $(x,y,c)$ with $(-x,-y,-c)$, we may assume that $x$ is positive. Likewise, since $y^p - x^p = c$ and $c\neq0$ (so $x\neq y$), we see that the absolute value of $c$ is minimized when $x$ and $y$ are consecutive. In particular, if $c>0$, then 
\[|c| = y^p - x^p \geq (x+1)^p - x^p = \sum_{k=0}^p x^k\binom pk - x^p \geq px^{p-1}.\]
Similarly, if $c < 0$, then 
\[|c| = x^p - y^p \geq x^p - (x-1)^p = \sum_{k=0}^p (x-1)^k\binom pk - (x-1)^p \geq p(x-1)^{p-1}.\]
Combining these cases implies that $|c| \geq p(|x| - 1)^{p-1}$ as claimed. \vspace{.2cm} 
\end{proof} 
Now we return to the proof of the Theorem \ref{thm:pthpower}.  
\begin{proof} 
Suppose that $\phi(x)=x^p+c$, that $\phi$ is irreducible, and that $\phi^3(a)=y^p$ for some $a,y\in\mathbb{Z}$. In particular, by replacing $(a,c)$ with $(-a, -c)$, we may assume that $a\geq 0$. From here, we proceed in cases: \vspace{.2cm} 

\textbf{Case (1):} Suppose first $c>0$ and $a\geq0$. Then $\phi(a) = a^p + c \geq c$ and so
\[\phi^2(a) = \phi(a)^p + c \geq c^p + c.\]
On the other hand, we have that $(\phi^2(a))^p +c = y^p$ and thus
\[c^p + c\leq\phi^2(a)=|\phi^2(a)|\leq \left(\frac{|c|}{p}\right)^{\frac{1}{p-1}} + 1\leq \left(\frac{|c|}{p}\right)^{\frac{1}{p-1}}+c\]
follows from Lemma \ref{lem:oddpbound}. However, this implies that $c=0$, and we reach a contradiction. 

In particular, we may assume that $c< 0$ and rewrite $\phi(x) = x^p - b$ for some $b> 0$. Moreover, since $b$ is not a perfect $p$th power (as $\phi$ is irreducible), it must be the case that $b$ is strictly in-between consecutive $p$th powers:
\[
\boxed{\phi(x)=x^p-b\;\;\text{and}\;\;\text{$s^p < b < (s+1)^p$}\;\;\text{for some $s>0$}} 
\]
From here, we make further cases.\vspace{.2cm} 
  
\textbf{Case (2):} Suppose $c<0$ and $0\leq a \leq s$. Then, since $\phi$ is an increasing function on the real line ($p$ is odd), we have that 
\[\phi^2(a)\leq \phi^2(s)=(s^p-b)^p-b<-b<0.\]
Therefore, $|\phi^2(a)|=-\phi^2(a)>b$. However, since $\phi(\phi^2(a))=y^p$, Lemma \ref{lem:oddpbound} implies 
\begin{equation} \label{eq:oddpcase2.2}
b<|\phi^2(a)|\leq \left(\frac{|c|}{p}\right)^{\frac{1}{p-1}} + 1=(b/p)^{1/(p-1)}+1,
\end{equation} 
Moreover, it is straightforward to check that the bound in \eqref{eq:oddpcase2.2} is violated for all $b\geq2$ and all $p\geq3$. Indeed, if $b\geq p$, then \eqref{eq:oddpcase2.2} implies that 
\[b< (b/p)^{1/(p-1)}+1\leq b/p+1\leq b/3+1\] 
and so $b\leq1$. Likewise, for $b< p$ the bound in \eqref{eq:oddpcase2.2} implies that 
\[b< (b/p)^{1/(p-1)}+1<2.\]
In particular, it follows that either $b=0$ or $b=1$, and we reach a contradiction (since both $x^p$ and $x^p+1$ are reducible polynomials). Therefore, it is not possible for $a\leq s$.\vspace{.2cm}    

\textbf{Case (3):} Next we consider the case when $c<0$ and $a \geq s + 2$. Then, since $\phi$ is an increasing function on the real line,   
\begin{equation}\label{eq:oddpcase3.1}
\phi(a) = a^p - c \geq (s+2)^p - (s+1)^p > p(s+1)^{p-1}.
\end{equation}
Moreover, it is clear that $p(s+1)^{p-1}>s+2$. In particular, we deduce that $\phi^2(a)>p(s+1)^{p-1}$, seen by replacing $a$ with $\phi(a)$ and repeating the argument above. On the other hand, since $\phi^2(a) = y^p$, we see from Lemma \ref{lem:oddpbound} that $\phi^2(a) < (b/p)^{1/(p-1)} + 1$. However, this can not possibly occur, since 
\begin{equation}\label{eq:oddpcase3.2}
p(s+1)^{p-1} \geq p\,b^\frac{p-1}{p} >\left(b/p\right)^\frac{1}{p-1} + 1.
\end{equation} 
Therefore, it is not possible that $\phi^3(a)=y^p$ when $c<0$ and $a\geq s+2$.\vspace{.2cm} 

In particular, combining the case above we see that $a=s+1$ and $c<0$. On the other hand, a similar argument obtained by replacing $a$ with $\phi(a)$ in cases 2 and 3 above shows that $\phi(a)=s+1$ also. For instance, if $\phi(a)=(s+1)^p-b\leq s$ (corresponding to Case 2 for $\phi(a)$ in place of $a$), then
\begin{equation}\label{oddp:phi(a)1}
\phi^2(a)\leq \phi(s)=s^p-b\leq s^p-(s+1)^p+s\leq-ps^{p-1}+s.
\end{equation} 
In particular, we see that 
\[|\phi^2(a)|=-\phi^2(a)\geq ps^{p-1}-s.\] 
However, as usual since $\phi^3(a)=y^p$, the bound in Lemma \ref{lem:oddpbound} implies that 
\[ ps^{p-1}-s\leq |\phi^2(a)|\leq (b/p)^{1/(p-1)}+1\leq \Big(\frac{(s+1)^p}{p}\Big)^{1/(p-1)}+1.\]
Moreover, it is straightforward to check that the bound above is violated unless $s=1$ and $p=3$. However, in this case we can instead use the lower bound $(s+1)^p-s^p-s\leq |\phi^2(a)|$ coming from \eqref{oddp:phi(a)1} to reach a contradiction. A similar adaptation of the argument in Case 3 obtained by replacing $a$ with $\phi(a)$ leads to a contradiction.        

In particular, we may deduce that $a=s+1=\phi(a)$ and $a$ is a fixed point of $\phi$. Therefore, $a=a^p-b$ and $-b=a-a^p$. On the other hand, since $a=\phi^3(a)=y^p$, we have that $a=y^p$ is a $p$th powered fixed point. Hence,  
\[\phi(x)=x^p-b=x^p+a-a^p=x^p+y^p-y^{p^2}\]
and $a=y^p$ are of the desired form (after a renaming of variables $y\rightarrow s$). 
\end{proof}

In particular, for odd $p$ we immediately obtain the following analog of Corollary \ref{cor:irre&non-special}. Namely, if $S$ is a set of polynomials of the form $x^p+c$ for $c\in\mathbb{Z}$, then $M_S$ contains many irreducible polynomials whenever $S$ contains an irreducible that is not of Type (I).  
\begin{cor}\label{cor:oddp+irre+nonspecial} Let $p\geq3$ be a prime and let $S=\{x^p+c_1,\dots,x^p+c_r\}$ for some $c_i\in\mathbb{Z}$. Moreover, assume that $S$ contains an irreducible polynomial $\phi$ that is not of Type (I). Then $\{\phi^3\circ f: f\in M_S\}$ is a set of irreducible polynomials. 
\end{cor}
\begin{proof} Combine Proposition \ref{prop:stability}, Theorem \ref{thm:pthpower}, and \cite[Corollary 30]{JonesIMRN}.  
\end{proof}
\section{The quadratic case}
In this section, we will construct an explicit subset (of positive proportion) in $M_S$ of irreducible polynomials when $S$ is a set of quadratic polynomials of the form $x^2+c$ for $c\in\mathbb{Z}$; see Theorem \ref{thm:main} from the introduction. In particular, we will assume throughout this section that at least one the polynomials in $S$ is irreducible in $\mathbb{Q}[x]$. Moreover, we will consider the case when $S$ contains at least two irreducible polynomials and the case when $S$ contains exactly one irreducible polynomial, separately.      

\subsection{Two irreducibles:} In this subsection, we will assume that $S$ contains at least two irreducible polynomials. In particular, we will strengthen the main result of \cite{Mathworks}. Note that by Corollary \ref{cor:irre&non-special} above, it suffices to assume that the given irreducible quadratics are either of Type (I) or Type (II), though not necessarily of the same type. With this in mind, we first classify the integral solutions to some related Diophantine equations.
\begin{prop}\label{prop:old1} Let $a,b,s,t\in\mathbb{Z}$ be such that 
\begin{equation}\label{eq:old1}
a^2+s^2-s^4=\pm{t^2}\;\;\;\;\text{and}\;\;\;\; b^2+t^2-t^4=\pm{s^2}.
\end{equation}
Moreover, assume that $s$ and $t$ are both non-zero. Then $s^2=t^2$. 
\end{prop}
\begin{proof} Assume that $s$ and $t$ are non-zero and that $(a,b,s,t)$ is an integral solution to \eqref{eq:old1} above. Moreover, suppose (for a contradiction) that $s^2>t^2$. Note that $s^2=t^2+1$ implies that $t=0$. In particular, $s^2>t^2+1$. Then
\begin{equation*}
(s^2-1)^2=s^4-s^2-(s^2-1)<s^4-s^2\pm t^2 =s^4-(s^2\pm t^2)<(s^2)^2.
\end{equation*}
In particular, $s^4-s^2+1 \pm t^2$ is an integer sandwiched between consecutive squares, and so cannot be a square itself. However, this contradicts the left side of \eqref{eq:old1}. Hence, it must be the case that $s^2\leq t^2$. Likewise, by symmetry, the right side of \eqref{eq:old1} implies that $t^2\leq s^2$, and we deduce that $s^2=t^2$ as claimed. 
\end{proof}
Similarly, we have the following result. 
\begin{prop}\label{prop:old2} Let $a,b,s,t\in\mathbb{Z}$ be such that 
\[a^2+s^2-s^4=\pm{(1+t^2)}\;\;\;\;\text{and}\;\;\;\; b^2-1-t^2-t^4=\pm{s^2}.\] 
Then $t=0$. 
\end{prop}
\begin{proof} First, we check that the systems given by $a^2+s^2-s^4=-(1+t^2)$ and $b^2-1-t^2-t^4=\pm{s^2}$ have no solutions $(a,b,s,t)$ modulo $4$.
%> m:=4; 
%> for x,y,s,t in Integers(m) do 
%> if x^2+s^2-s^4 eq -(1+t^2) and y^2-1-t^2-t^4 eq -(s^2) then 
%> print <s,t>; 
%> end if; 
%> end for;
In particular, it suffices to consider $a^2+s^2-s^4=1+t^2$ and $b^2-1-t^2-t^4=\pm{s^2}$. Let us assume first that 
\begin{equation}\label{eq:mixed1}
a^2+s^2-s^4=(1+t^2)\;\;\;\;\text{and}\;\;\;\; b^2-1-t^2-t^4=\pm s^2.
\end{equation}
Suppose for a contradiction that $s^2>t^2>0$. Note that if $s^2=t^2+1$, then $t=0$. Hence, we may assume that $s^2-t^2\geq2$. Then we see that
\begin{equation*}
(s^2-1)^2=s^4-s^2+1-s^2<s^4-s^2+1+t^2=s^4-(s^2-t^2-1)<(s^2)^2. 
\end{equation*} 
In particular, $s^4-s^2+1+t^2$ is sandwiched between consecutive squares and so cannot be a square itself. However, this contradicts the left side of \eqref{eq:mixed1}. Likewise if $t^2>s^2$,
\begin{align*}
(t^2+1)^2=t^4+t^2+1+t^2>t^4+t^2+1\pm s^2>(t^2)^2. 
\end{align*} 
In particular, $t^4+t^2+1\pm s^2$ is sandwiched between consecutive squares, contradicting the right side of \eqref{eq:mixed1}. Therefore, we deduce that either $t=0$ or $s^2=t^2$. However, in the latter case the left side of \eqref{eq:mixed1} implies that $(x-t^2)(x+t^2)=1$, from which we again deduce that $t=0$ as desired. 
\end{proof}
Likewise, we have the following result.
\begin{prop}\label{prop:old3} Let $a,b,s,t\in\mathbb{Z}$ be such that 
\begin{equation}\label{eq:old3}
  a^2-1-s^2-s^4=\pm{(1+t^2)}\;\;\;\;\text{and}\;\;\;\; b^2-1-t^2-t^4=\pm{(1+s^2)}.  
\end{equation}
Moreover, assume that $s$ and $t$ are both non-zero. Then $s^2=t^2$. 
\end{prop}
\begin{proof} 
Assume that $s$ and $t$ are non-zero and that $(a,b,s,t)$ is an integral solution to \eqref{eq:old3} above. Moreover, suppose (for a contradiction) that $s^2>t^2$. Note first that $s^2=t^2+1$ implies that $t=0$. Likewise, $s^2=t^2+2$ has no solutions modulo $4$. Therefore, we may assume that $s^2>t^2+2$. Hence, 
\begin{align*}
(s^2+1)^2&=s^4+s^2+1+s^2>s^4+s^2+1+(t^2+1)\\[3pt]
&\geq s^4+s^2+1\pm(t^2+1)=s^4+(s^2+1\pm(t^2+1))\\[3pt]  
&>s^4+(t^2+2\pm(t^2+1))>(s^2)^2. 
\end{align*} 
In particular, $s^4+s^2+1\pm(t^2+1)$ is sandwiched between consecutive squares, contradicting the left side of \eqref{eq:mixed1}. Therefore, $t^2\geq s^2$. However, by symmetry the right side of \eqref{eq:old3} implies that $s^2\geq t^2$. Hence, $s^2=t^2$ as claimed. 
\end{proof}

We now have the tools in place to prove the following irreducibility result for sets $S$ containing at least 2 irreducible polynomials, both of which are either of Type (I) or of Type (II).  

\begin{thm}\label{thm:2irreducble,special}
Let $S=\{x^2+c_1,\dots,x^2+c_r\}$ for some distinct $c_1,\dots,c_r\in\mathbb{Z}$ and let $\phi_i(x)=x^2+c_i$. Moreover, assume that $\phi_1$ and $\phi_2$ are both irreducible in $\mathbb{Q}[x]$ and both are of Type (I) or of Type (II), though not necessarily of the same type. Then either 
\begin{equation}\label{eq:bothspecial}
\{\phi_1^3\circ\phi_2\circ f\, :\, f\in M_S\}
\quad
\text{or}
\quad
\{\phi_2^3\circ\phi_1\circ f\, :\, f\in M_S\}
\end{equation}
is a set of irreducible polynomials in $\mathbb{Q}[x]$.    
\end{thm}
\begin{proof} Assume that $\phi_1$ and $\phi_2$ are both irreducible in $\mathbb{Q}[x]$ and both are of Type (I) or of Type (II). In particular, $\phi_1^3$ and $\phi_2^3$ are both irreducible by \cite[Corollary 1.3]{Stoll}. Hence, if a polynomial of the form $\phi_1^3\circ\phi_2\circ f$ is reducible for some $f\in M_S$ (including the case when $f$ is the identity) \emph{and} a polynomial of the form $\phi_2^3\circ\phi_1\circ g$ is reducible for some $g\in M_S$ (including the case when $g$ is the identity), then Proposition \ref{prop:stability}, Lemma \ref{lem:family1refinement}, and Lemma \ref{lem:family2refinement} together imply that we obtain an integral solution $(a,b,s,t)$ to one of the systems of equations in Propositions \ref{prop:old1}, \ref{prop:old2}, or \ref{prop:old3}. However, a contradiction occurs in each of the three cases: Proposition \ref{prop:old1} implies that $\phi_1=\phi_2$, Proposition \ref{prop:old2} implies that either $\phi_1$ or $\phi_2$ is equal to $x^2-1$ (and so reducible), and Proposition \ref{prop:old1} implies that $\phi_1=\phi_2$. Therefore, one of the sets in \eqref{eq:bothspecial} must consist entirely of irreducible polynomials as claimed.         
\end{proof}
In particular, by combining Corollary \ref{cor:irre&non-special} and Theorem \ref{thm:2irreducble,special}, we deduce the following improvement of the main result in \cite{Mathworks}. 
\begin{cor} Let $S=\{x^2+c_1,\dots,x^2+c_r\}$ for some distinct $c_1,\dots,c_r\in\mathbb{Z}$ and let $\phi_i(x)=x^2+c_i$. Moreover, assume that at least two of the elements of $S$, say $\phi_1$ and $\phi_2$, are irreducible in $\mathbb{Q}[x]$. Then at least one of the sets: \vspace{.1cm}  
\[
\begin{array}{lr}
\{\phi_1^4\circ f\, :\, f\in M_S\},& \{\phi_2^4\circ f\, :\, f\in M_S\},\\[5pt]
\{\phi_1^3\circ\phi_2\circ f\, :\, f\in M_S\},& \{\phi_2^3\circ\phi_1\circ f\, :\, f\in M_S\}
\end{array}
\vspace{.075cm}
\]
is a set of irreducible polynomials in $\mathbb{Q}[x]$.    
\end{cor}
\begin{proof} Suppose that $S$ contains two irreducible polynomials $\phi_1$ and $\phi_2$. If either $\phi_1$ or $\phi_2$ is is not of Type (I) or of Type (II) (say $\phi_1$ is not), then Corollary \ref{cor:irre&non-special} implies that $\{\phi_1^4\circ f\,:\, f\in M_S\}$ is a set of irreducible polynomials. Therefore, we can assume that both $\phi_1$ and $\phi_2$ are of Type (I) or of Type (II). In this case, Theorem \ref{thm:2irreducble,special} implies that $\{\phi_i^3\circ\phi_j\circ f\, :\, f\in M_S\}$ is a set of irreducible polynomials for $(i,j)=(1,2)$ or $(i,j)=(2,1)$ as claimed. 
\end{proof}

\subsection{One irreducible of Type (I)}
As for the case when $S$ contains exactly only irreducible polynomial, it suffices to assume that this polynomial is either of Type (I) or of Type (II) by Corollary \ref{cor:irre&non-special}. In particular, we study the Type I case in this subsection and the Type II case in the following subsection. Moreover, since the case of iterating a single \emph{integral} quadratic $x^2+c$ is already well understood \cite{Stoll}, we may assume that $S$ contains an additional reducible polynomial, which is necessarily of the form $x^2-t^2$ for some $t\in\mathbb{Z}$. In particular, we will assume that $S$ contains polynomials of the form $x^2+s^2-s^4$ and $x^2-t^2$ for some $s,t\in\mathbb{Z}$ throughout this subsection. With this in mind, we first classify the integral solutions to a related Diophantine equation. 
\begin{prop}\label{prop:case1} Let $y,s,t\in\mathbb{Z}$ be such that 
\begin{equation}\label{eq:case1}
(y^2+s^2-s^4)^2-t^2=\pm{s}^2.
\end{equation}
Then $s=0,\pm{1}$. 
\end{prop} 
\begin{proof} Let $(y,s,t)$ be a solution to one of the equations in \eqref{eq:case1} above and assume that $s\neq0,\pm{1}$. Moreover, note that by replacing $y$, $s$, or $t$ with its negative, we may assume that $y,s,t\geq0$. Next, we factor and take absolute values of \eqref{eq:case1}
\[|A_1|\cdot |A_2|:=|(y^2+s^2-s^4-t)|\cdot|(y^2+s^2-s^4+t)|=s^2\]
to deduce that $\max\{|A_1|, |A_2|\}\leq s^2$. Therefore, 
\begin{equation}\label{ineq1}
2|y^2+s^2-s^4|\leq|y^2+s^2-s^4-t|+|y^2+s^2-s^4+t|\leq 2 s^2
\end{equation}
follows from the triangle inequality. Likewise, we then see that 
\[|y^2-s^4|-|s^2|=|y^2-s^4|-|-s^2|\leq|y^2+s^2-s^4|\leq s^2\]
follows from the fact that $|a-b|\geq |a|-|b|$ and the inequality above. In particular, we deduce that 
\begin{equation}\label{ineq2}
|y-s^2|\cdot(y+s^2)=|y-s^2|\cdot|y+s^2|=|y^2-s^4|\leq 2s^2,
\end{equation}
since $y+s^2$ is positive as $y\geq0$. From here, we proceed in cases, obtaining a contradiction in each.\\[5pt]
\textbf{Case(1):} Suppose that $|y-s^2|\geq2$. Then \eqref{ineq2} implies that $2s^2+2y\leq 2s^2$ so that $y=0$. However, in this case \eqref{ineq1} implies that $|1-s^2|\leq 1$, contradicting our assumption that $s\geq2$. \\[5pt]
\textbf{Case(2):} Suppose that $|y-s^2|=1$. However, if $y-s^2=1$, then \eqref{ineq2} implies that $2s^2+1\leq 2s^2$, a contradiction. Therefore, we may assume that $y=s^2-1$ in this case. In particular, \eqref{eq:case1} implies that $s^4 - 2s^2 - t^2 + 1=\pm{s^2}$. We handle these subcases sequentially. 

First, if $s^4 - 2s^2 - t^2 + 1=s^2$, then we obtain a rational point $(s,t)$ on the hyperelliptic curve $C_1: t^2=s^4 - 3s^2 + 1$. However, we compute with {\tt{Magma}} that $C_1$ is birational over $\mathbb{Q}$ with the elliptic curve $E_1$ in Weierstrass form $E_1: Y^2 = X^3 + 6X^2 + 5X$. Moreover, we compute that $E_1(\mathbb{Q})\cong\mathbb{Z}/2\mathbb{Z}\times \mathbb{Z}/2\mathbb{Z}$, from which it follows that the only affine rational points on $C_1$ correspond to those with $s=0$, a contradiction. 

Similarly, if $s^4 - 2s^2 - t^2 + 1=-s^2$, then we obtain a rational point $(s,t)$ on the hyperelliptic curve $C_2: t^2=s^4 - s^2 + 1$. However, we compute with {\tt{Magma}} that $C_2$ is birational over $\mathbb{Q}$ with the elliptic curve $E_2$ in Weierstrass form $E_2: Y^2 = X^3 + 2X^2 - 3X$. Moreover, we compute that $E_2(\mathbb{Q})\cong\mathbb{Z}/2\mathbb{Z}\times \mathbb{Z}/4\mathbb{Z}$, from which it follows that the only affine rational points on $C_2$ correspond to those with $s\in\{0,\pm{1}\}$, a contradiction. 
\\[5pt]
\textbf{Case(3):} Suppose that $y=s^2$. Then the original equation \eqref{eq:case1} implies that $s^4-t^2=\pm{s}^2$. However, when $s^4-t^2=\pm{s}^2$, simply rearranging and factoring yields $s^2(s^2\pm 1)=t^2$. In particular, since $s\neq0$ it must be the case that $s^2\pm 1$ is itself a square: $s^2\pm 1=z^2$ for some integer $z$. However, again factoring this new equation, $(s-z)(s+z)=s^2-z^2=\pm1$, we easily deduce that $s=0,\pm1$ and so reach a contradiction. 
\end{proof}

In particular, we have the tools in place to prove the following irreducibility result for sets $S$ containing an irreducible polynomial of Type (I) and a second reducible polynomial (meeting a further stipulation).

\begin{thm}\label{thm:family1+reducible+stipulation}
Let $S=\{x^2+c_1,\dots,x^2+c_r\}$ for some distinct $c_1,\dots,c_r\in\mathbb{Z}$ and let $\phi_i(x)=x^2+c_i$. Moreover, assume that $\phi_1$ is irreducible and of Type (I), that $\phi_2$ is reducible, and write $c_1=s^2-s^4$ and $c_2=-t^2$ for some $s,t\in\mathbb{Z}$. If $s^2\neq t^2$, then  
\[
\{\phi_1^3\circ\phi_2\circ\phi_1\circ f\, :\, f\in M_S\}
\]
is a set of irreducible polynomials in $\mathbb{Q}[x]$. 
\end{thm}
\begin{proof} Assume that $c_1$ and $c_2$ are as above, that $\phi_1$ is irreducible (so $s\neq0,\pm{1}$), and that $t^2\neq s^2$. Then $\phi_1^3$ is irreducible by \cite[Corollary 1.3]{Stoll}. Next, if $\phi_1^3(\phi_2(0))$ is a square, then $-t^2=\phi_2(0)=\pm{s^2}$ by Lemma \ref{lem:family1refinement}. However, clearly -$t^2=s^2$ is impossible and $-t^2\neq -s^2$ by assumption. Hence, $\phi_1^3\circ\phi_2$ is irreducible by Proposition \ref{prop:stability}. Likewise, if a polynomial of the form $\phi_1^3\circ\phi_2\circ\phi_1\circ f$ for some $f\in M_S$ is reducible, then Proposition \ref{prop:stability} implies that $\phi_1^3\circ\phi_2\circ\phi_1(a)$ is a square for some $a\in\mathbb{Z}$. But then Lemma \ref{lem:family1refinement} implies that $\phi_2(\phi_1(a))=\pm{s^2}$ and Proposition \ref{prop:case1} implies that $s=0,\pm{1}$, a contradiction. Therefore, every polynomial of the form $\phi_1^3\circ\phi_2\circ\phi_1\circ f$ for $f\in M_S$ is irreducible as claimed. 
\end{proof}
So it suffices to consider the case when $S$ contains polynomials of the form $x^2+t^2-t^4$ and $x^2-t^2$ for some $t\in\mathbb{Z}$ (i.e., $s^2=t^2$) in this subsection. To construct a positive proportion of irreducible polynomials in this case, we begin with the following result. 
\begin{prop}\label{prop:irreoctic}
Let $\phi_1(x)=x^2+t^2-t^4$ and $\phi_2(x)=x^2-t^2$ for some $t\in\mathbb{Z}\setminus\{0,\pm{1}\}$. Then $\phi_1\circ\phi_1\circ\phi_2$ is irreducible over $\mathbb{Q}$.
\end{prop}
We break this proof into two pieces. First, we note that $\phi_1\circ\phi_1\circ\phi_2$ is of the form $f(x^2)$, where $f(x)=x^4-4t^2x^3+(4t^4+2t^2)x^2-4t^4x+t^2$. In particular, we first show that such quartics $f$ are irreducible. 
\begin{lem}\label{lem:quartic} Let $f(x)=x^4-4t^2x^3+(4t^4+2t^2)x^2-4t^4x+t^2$ for some $t\in\mathbb{Z}\setminus\{0,\pm{1}\}$. Then $f$ is irreducible over $\mathbb{Q}$. 
\end{lem}
\begin{proof}
Assume for the sake of contradiction that $f(x)$ is reducible. We first compute the reduced form of $f$; that is, an equivalent polynomial to $f$ that has no cubic term. We note that 
\[f(x+t^2) = x^4 + (2t^2 - 2t^4)x^2 + (t^8 - 2t^6 + t^2). \]
Note that clearly, $f$ is irreducible if and only if $f(x + t^2)$ is irreducible. We now note that by \cite[Theorem 2]{brookfield2007factoring}, the quartic $f(x+t^2)$ is reducible if and only if either $(2t^2 - 2t^4)^2 - 4(t^8 - 2t^6 + t^2)$ is a perfect square, or if $t^8 - 2t^6 + t^2$ is a perfect square. We first note that
\[(2t^2 - 2t^4)^2 - 4(t^8 - 2t^6 + t^2) = 4t^4 - 4t^2, \]
which is one less than $(2t^2 - 1)^2 = 4t^4 - 4t^2 + 1$. Thus, if it were a perfect square, it must be 0, which implies that $t = 0, \pm 1$, a contradiction. Now, note that 
\begin{align*}
    (t^4 - t - 1)^2 &= t^8 - 2t^5 + t^2 - 2t + 1 > t^8-2t^6+t^2 > (t^4-t-2)^2
\end{align*}
for all $|t| \geq 3$. We can quickly confirm that $t = \pm 2$ does not yield a perfect square. Thus $t^8 - 2t^6 + t^2$ cannot be a square for the given domain. Combining the two cases yields that $f(x+t^2)$ is irreducible, and therefore, $f(x)$ must be irreducible, as desired. 
\end{proof}
Next, we show that in this special case, $f(x)$ is irreducible if and only if $f(x^2)$ is irreducible. This involves finding all the rational points on a certain hyperelliptic curve of genus two. 
\begin{lem}\label{lem:octic}
Let $f(x)=x^4-4t^2x^3+(4t^4+2t^2)x^2-4t^4x+t^2$ for some $t\in\mathbb{Q}\setminus\{0,\pm{1}\}$. Then $f(x)$ is irreducible in $\mathbb{Q}[x]$ if and only if $f(x^2)$ is irreducible in $\mathbb{Q}[x]$. 
\end{lem}
\begin{proof} Clearly if $f(x)$ is reducible then so is $f(x^2)$. On the other hand, if $f(x)$ is irreducible and $f(x^2)$ is reducible, then \cite[Proposition 2.4]{Octics} implies that there is an $x\in\mathbb{Q}$ such that 
\[x^4 + (-8t^4 - 4t^2 - 12t)x^2 + (-32t^4 - 32t^3)x + 16t^8 - 48t^6 -
    16t^5 + 4t^4 - 8t^3 + 4t^2=0.\]
Moreover, a point search yields the following four rational solutions 
\[(x,t)\in\{(0,0), (0,-1), (-2,1), (6,1)\}.\] 
On the other hand, we compute with {\tt{Magma}} \cite{Magma} that the equation above is birational over $\mathbb{Q}$ with the hyerelliptic curve
\[C: Y^2=2X^5 + 6X^4 - 4X^3 - 12X^2 + 8\]
of genus $2$. Moreover, we compute that the rational points on the Jacobian of $C$ have rank one, and so the method of Chabauty, together with the Mordell-Weil sieve, can be used to determine the full set of rational point of $C$. In fact, for genus $2$ hyperelliptic curves this procedure has been implemented in {\tt{Magma}}, and we compute that $\# C(\mathbb{Q})=4$ in this case. In particular, the only rational values of $t$ for which $f(x^2)$ can be reducible (while $f$ is irreducible) are $t=0,\pm{1}$. The claim follows; see \href{https://sites.google.com/a/alumni.brown.edu/whindes/research/code}{Mathworks 2023} for the relevant {\tt{Magma}} \cite{Magma} code.  
%> R<x,t>:=PolynomialRing(Rationals(),2); 
%> a:=-4*t^2; 
%> b:=4*t^4+2*t^2;
%> c:=-4*t^4; 
%> n:=t; 
%> G:=x^4-(2*b+12*n)*x^2+(8*c+8*a*n)*x+b^2-4*a*c-4*b*n+4*n^2; 
%> A<x,t>:=AffineSpace(Rationals(),2); 
%> C:=Curve(A,G); 
%> Genus(C); 
%> T,C,p:=IsHyperelliptic(C); 
%> C; 
%> J:=Jacobian(C); 
%> RankBound(J); 
%> Pts:=Points(C: Bound:=1000); 
%> P:=J![Pts[1],Pts[3]]; 
%> Order(P); 
%> Chabauty(P); /// so we indeed have found all of the points of C
%> for x, t in [-100..100] do 
%> a:=-4*t^2; 
%> b:=4*t^4+2*t^2;
%> c:=-4*t^4; 
%> n:=t; 
%> G:=x^4-(2*b+12*n)*x^2+(8*c+8*a*n)*x+b^2-4*a*c-4*b*n+4*n^2; 
%> if G eq 0 then 
%> print <x,t>; 
%> end if; 
%> end for; 
%// so we’ve found four rational points on the original curve too, and so this must be all of them. 
\end{proof}
Therefore, we see that $\phi_1^2\circ\phi_2$ is irreducible by combining the previous two lemmas. Next, we classify the integral solutions to a related Diophantine equation. 
\begin{lem}\label{lem:lastpairnosquares}
Let $\phi_1(x)=x^2+t^2-t^4$ and $\phi_2(x)=x^2-t^2$ for some $t\in\mathbb{Z}\setminus\{0,\pm{1}\}$. Then $\phi_1^2\circ\phi_2\circ\phi_1(a)\neq y^2$ for all $a,y\in \mathbb{Z}$.       
\end{lem}
\begin{proof}
Suppose for a contradiction that $\phi_1^2\circ\phi_2\circ\phi_1(a)= y^2$ for some $a,y\in \mathbb{Z}$ and let $b=\phi_1\circ\phi_2\circ\phi_1(a)$. Then $\phi_1(b)=y^2$ so that $(b-y)(b+y)=t^4-t^2$. In particular, $|b|\leq |t^4-t^2|=t^4-t^2$ follows from the triangle inequality; here we use that $t^4-t^2\neq0$ since $t\neq 0,\pm{1}$. Note that by replacing $a$ and $t$ with their negatives, we may assume that $a,t\geq0$. From here we proceed in cases. Suppose first that $a\geq t^2+1$. Then $\phi_1(a)\geq \phi_1(t^2+1)=3t^2+1>0$ and so
\[\phi_2\circ\phi_1(a)\geq \phi_2(3t^2+1)>\phi_2(3t^2)=9t^4-t^2> t^4.\]
In particular, we deduce that 
\[b=\phi_1\circ\phi_2\circ\phi_2(a)>\phi_1(t^4)=t^{16}-t^4+t^2>t^4-t^2\]
for all $t\in\mathbb{Z}\setminus\{0,\pm{1}\}$, a contradiction. Likewise when $a \leq t^2-2$, we see 
\[\phi_1(a)\leq\phi(t^2-2)=-3t^2+4<0.\]
Hence $-\phi_1(a)>3t^2-4>0$. But then 
\[\phi_2\circ\phi_1(a)=\phi_2(-\phi_1(a))>\phi_2(3t^2-4)=9t^4-25t^2+16>t^4>0\]
for all $t\in\mathbb{Z}\setminus\{0,\pm{1}\}$. Therefore, $\phi_1(b)>t^8-t^4+t^2>t^4-t^2$, and we again reach a contradiction. Hence, it remains to consider $a=t^2,t^2-1$. However, we compute that 
\[\phi_1\circ\phi_2\circ\phi_1(t^2)=t^8 - 2t^6 + t^2>t^4-t^2\]
for all $t\in\mathbb{Z}\setminus\{0,\pm{1}\}$. Hence, $a\neq t^2$ also. Lastly, when $a=t^2-1$ we see that 
\[b=\phi_1\circ\phi_2\circ\phi_1(t^2-1)=t^8 - 6t^6 + 10t^4 - 5t^2 + 1>t^4-t^2\]
for all $t\in\mathbb{Z}\setminus\{0,\pm{1}\}$, a contradiction. The result follows.             
\end{proof}
In particular, we conclude this subsection with the following consequence. 
\begin{thm}\label{thm:family1+s=t}
Let $S=\{x^2+c_1,\dots,x^2+c_r\}$ for some distinct $c_1,\dots,c_r\in\mathbb{Z}$ and let $\phi_i(x)=x^2+c_i$. Moreover, assume that $c_1=t^2-t^4$ and $c_2=-t^2$ for some $t\in\mathbb{Z}$ and that $\phi_1$ is irreducible in $\mathbb{Q}[x]$. Then  
\[
\{\phi_1^2\circ\phi_2\circ\phi_1\circ f\, :\, f\in M_S\}
\]
is a set of irreducible polynomials in $\mathbb{Q}[x]$. 
\end{thm}
\begin{proof} First, Proposition \ref{prop:irreoctic} implies that $\phi_1^2\circ\phi_2$ is irreducible. Hence, if a polynomial of the form $\phi_1^2\circ\phi_2\circ\phi_1\circ f$ for some $f\in M_S$ (or $f$ is the identity) is reducible, then it follows from Proposition \ref{prop:stability} that $\phi_1^2\circ\phi_2\circ\phi_1(a)=y^2$ for some $a,y\in\mathbb{Z}$. But then Lemma \ref{lem:lastpairnosquares} implies that $t=0,\pm{1}$ and $c_1=0$, contradicting our assumption that $\phi_1$ is irreducible. 
\end{proof}
\subsection{One irreducible of Type (II)} In this subsection, we will assume that $S$ contains polynomials of the form $x^2-1-s^2-s^4$ and $x^2-t^2$ for some $s,t\in\mathbb{Z}$. With this in mind, we first classify the integral solutions to some related Diophanitne equations.
\begin{prop}\label{prop:case2} Let $y,s,t\in\mathbb{Z}$ be such that
\begin{equation}\label{eq:case2}
(y^2-1-s^2-s^4)^2-t^2=\pm(1+s^2).
\end{equation}
Then $s=0$.  
\end{prop} 
\begin{proof} Let $(y,s,t)$ be a solution to one of the equations in \eqref{eq:case1} and assume that $s\neq0$. Note that if $s=\pm{1}$, then we obtain an integral solution to $a^2-t^2=\pm{2}$, which has no solutions modulo $4$. Therefore, we may assume that $s\neq\pm{1}$ also. Likewise, by replacing $y$, $s$, or $t$ with its negative, we may assume that $y,s,t\geq0$. Next, we factor and take absolute values of \eqref{eq:case2}
\[|B_1|\cdot |B_2|:=|(y^2-1-s^2-s^4-t)|\cdot|(y^2-1-s^2-s^4+t)|=s^2+1\]
to deduce that $\max\{|B_1|, |B_2|\}\leq s^2+1$. In particular, it follows that from the triangle inequality that 
\[
|y^2-1-s^2-s^4|\leq s^2+1.
\]
On the other hand, since $|a-b|\geq|a|-|b|$, we have that 
\[|y^2-s^4|-|s^2+1|\leq|(y^2-s^4)-(s^2+1)|\leq s^2+1.\]
Finally, by factoring (and using that $y\geq0$) we deduce that
\begin{equation}\label{ineq3}
|y-s^2|\cdot (y+s^2) \leq 2s^2+2.  
\end{equation}
From here, we proceed in cases, obtaining a contradiction in each.\\[5pt]
\textbf{Case(1):} Suppose that $|y-s^2|\geq3$. Then \eqref{ineq3} implies that $3y+3s^2\leq 2s^2+2$. In particular, $3y+s^2\leq 2$ so that $s\in\{0,\pm{1}\}$, a contradiction.  \\[5pt]
\textbf{Case(2):} Suppose that $|y-s^2|=2$. If $y-s^2=2$, then \eqref{ineq3} implies that $4s^2+4\leq 2s^2+2$, and we obtain a contradiction. Likewise, if $y-s^2=-2$, then \eqref{ineq3} implies that $4s^2-4\leq 2s^2+2$. However, this implies that $s^2\leq3$, and hence $s\in\{0,\pm{1}\}$, a contradiction.
\\[5pt]
\textbf{Case(3):} Suppose that $|y-s^2|=1$. If $y-s^2=1$, then \eqref{eq:case2} implies that $s^4-t^2=\pm{(s^2+1)}$. However, the equation $s^4-t^2=s^2+1$ has no solutions modulo $4$. Finally, when $s^4-t^2=-(s^2+1)$ then we obtain a rational point on the hyperelliptic curve $C_3: t^2=s^4+s^2+1$. Moreover, we compute with {\tt{Magma}} that $C_3$ is birational over $\mathbb{Q}$ to the elliptic curve $E_3$ with Weierstrass equation $E_3: Y^2=X^3-2X^2-3X$. Furthermore, we compute that $E_3(\mathbb{Q})\cong \mathbb{Z}/2\mathbb{Z}\times\mathbb{Z}/2\mathbb{Z}$, from which it follows that the only affine rational points on $C_3$ are those corresponding to $s=0$, a contradiction.

Likewise, if $y-s^2 = -1$, then \eqref{eq:case2} implies that $9s^4-t^2=\pm{(s^2+1)}$. However, $9s^4-t^2 = s^2+1$ has no solutions modulo 4. Therefore, it suffices to rule out $9s^4-t^2 = -(s^2+1)$. In this case, we obtain a rational point on the hyperelliptic curve $C_4: t^2=9s^4-s^2-1$. Moreover, we compute with {\tt{Magma}} that $C_4$ is birational over $\mathbb{Q}$ to the elliptic curve $E_4$ with Weierstrass equation $E_4: Y^2=X^3+2/9X^2+37/81X$. Furthermore, we compute that $E_4(\mathbb{Q})\cong \mathbb{Z}/2\mathbb{Z}$, from which it follows that there are no affine points in $C_4(\mathbb{Q})$, a contradiction.  
\\[5pt] 
\textbf{Case(4):} Suppose that $y=s^2$. Then the original equation \eqref{eq:case2} implies that $s^4 + 2s^2 - t^2 + 1=\pm{(1+s^2)}$. However, when $s^4 + 2s^2 - t^2 + 1=1+s^2$ we see that $s^2(s^2+1)=t^2$. In particular, since $s\neq0$ it must be the case that $s^2+1=z^2$ for some $z\in\mathbb{Z}$. However, by factoring $(s-z)(s+z)=1$ one easily deduces that $(s,z)\in\{(0,\pm{1})\}$, a contradiction. Likewise, the equation $s^4 + 2s^2 - t^2 + 1=-(1+s^2)$ has no solutions modulo 4. 
\end{proof}

In particular, we have the tools in place to prove the following irreducibility result for sets $S$ containing an irreducible polynomial of Type (II) and a second reducible polynomial.

\begin{thm}\label{thm:family2+reducible}
Let $S=\{x^2+c_1,\dots,x^2+c_r\}$ for some distinct $c_i\in\mathbb{Z}$ and let $\phi_i(x)=x^2+c_i$. Moreover, assume that $\phi_1$ is of Type(II) and irreducible and that $\phi_2$ is reducible. Then  
\[
\{\phi_1^3\circ\phi_2\circ\phi_1\circ f\, :\, f\in M_S\}
\]
is a set of irreducible polynomials in $\mathbb{Q}[x]$. 
\end{thm}
\begin{proof} Write $c_1=-1-s^2-s^4$ and $c_2=-t^2$ for some $s,t\in\mathbb{Z}$ and assume that $\phi_1$ is irreducible (so $s\neq0$). Then $\phi_1^3$ is irreducible by \cite[Corollary 1.3]{Stoll}. Next, if $\phi_1^3(\phi_2(0))$ is a square, then $-t^2=\phi_2(0)=\pm{(s^2+1)}$ by Lemma \ref{lem:family2refinement}. However, clearly -$t^2=s^2+1$ is impossible and $-t^2=-(s^2+1)$ implies that $s=0$, a contradiction. Therefore, $\phi_1^3\circ\phi_2$ is irreducible by Proposition \ref{prop:stability}. Likewise, if a polynomial of the form $\phi_1^3\circ\phi_2\circ\phi_1\circ f$ for some $f\in M_S$ is reducible, then Proposition \ref{prop:stability} implies that $\phi_1^3\circ\phi_2\circ\phi_1(a)$ is a square for some $a\in\mathbb{Z}$. But then Lemma \ref{lem:family2refinement} implies that $\phi_2(\phi_1(a))=\pm{(s^2+1)}$ and Proposition \ref{prop:case2} implies that $s=0$, a contradiction. Therefore, every polynomial of the form $\phi_1^3\circ\phi_2\circ\phi_1\circ f$ for $f\in M_S$ is irreducible as claimed.              
\end{proof} 
 
\subsection{Main quadratic result}
We now have all of the pieces in place to prove our main result in the quadratic case: when $S$ is a set of quadratic polynomials of the form $x^2+c$ and the $c$'s are integral, the semigroup $M_S$ generated by $S$ contains a positive proportion of irreducible polynomials (over $\mathbb{Q}$) if and only if it contains at least one irreducible polynomial. However, before we do this we first prove that the semigroups we consider in this paper are free; see \cite[Theorem 3.1]{DiscCont} for a similar statement. 
\begin{lem}\label{lem:free} Let $p$ be a prime, let $K$ be a field of characteristic not equal to $p$, and let $S=\{x^p+c_1, \dots, x^p+c_r\}$ for some $c_i\in K$. Then the semigroup $M_S$ is free, i.e., if 
\begin{equation}\label{eq:free}
\theta_1\circ\dots\circ \theta_n=\tau_1\circ\dots\circ \tau_m
\end{equation}
for some $\theta_i,\tau_j\in S$ and $n,m\geq1$, then $n=m$ and $\theta_i=\tau_i$ for all $1\leq i\leq n$. 
\end{lem}
\begin{proof} Note first that if \eqref{eq:free}, then equating degrees implies that $p^n=p^m$ and so $m=n$. The claim above is clear when $n=1$. On the other hand, when $n>1$ write $\theta_1=x^p+c_1$, $\tau_1=x^p+c_1'$, $A=\theta_2\circ\dots\circ \theta_n$, and $B=\tau_2\circ\dots\circ \tau_m$. Then \eqref{eq:free} implies that $A^p+c_1=B^p+c_1'$. Hence, 
\[A^p-B^p=(A-B)(A^{p-1}+A^{p-2}B+\dots+AB^{p-2}+B^{p-1})=c_1-c_1'.\] 
However, since $A$ and $B$ are monic polynomials of equal degree $p^{n-1}$, we see that the leading coefficient of $G_{A,B}:=A^{p-1}+A^{p-2}B+\dots+AB^{p-2}+B^{p-1}$ is $p$. In particular, since $K$ has characteristic distinct form $p$, it must be the case that $G_{A,B}$ is non-constant. Therefore, $c_1=c_1'$, since otherwise $G_{A,B}$ is a unit in $K[x]$. Thus $\theta_1=\tau_1$ and $A-B=0$ by above. Therefore, $\theta_2\circ\dots\circ \theta_n=A=B=\tau_2\circ\dots\circ \tau_m$ also. In particular, repeating the argument above with the remaining strings yields the desired conclusion. 
\end{proof}
\begin{rem}\label{rem:length} In what follows, if $f=\theta_1\circ\dots\circ \theta_n$ for some $\theta_i\in S$ is an arbitrary element of $M_S$, then we define the \emph{length} of $f$ to be $|f|=n$, which is well-defined by Lemma \ref{lem:free}. In particular, it is clear that $f\in M_S$ satisfies $\deg(f)\leq B$ if and only if $|f|\leq \log_2(B)$. Hence, we may replace $\deg(f)$ with $|f|$ in Definition \ref{def:proportion} without effect. 
\end{rem}
\begin{proof}[(Proof of Theorem \ref{thm:main})] Statements (1)-(6) are restatements of the following results in order: \cite[Corollary 1.3]{Stoll}, Corollary \ref{cor:irre&non-special}, Theorem \ref{thm:2irreducble,special}, Theorem \ref{thm:family1+reducible+stipulation}, Theorem \ref{thm:family1+s=t}, and Theorem \ref{thm:family2+reducible}. In particular, in all cases we have constructed a polynomial $g\in M_S$ with $|g|\leq 5$ such that $\{g\circ f: f\in M_S\}$ is a set of irreducible polynomials. 
In particular \vspace{.15cm}
\[\liminf_{B\rightarrow\infty}\,\frac{\#\big\{f\in M_S\,: \deg(f)\leq B\;\text{and $f$ is irreducible}\big\}}{\#\big\{f\in M_S\,: \deg(f)\leq B\big\}}\geq \frac{1}{(\#S)^5}>0 
\vspace{.15cm}
\]
 follows from Lemma \ref{lem:free} and Remark \ref{rem:length} above. 
\end{proof}
\section{The odd prime case}
In this section, we will assume that $p$ is an odd prime and that $S$ is a set of polynomials of the form $x^p+c$ for some $c\in\mathbb{Z}$. Moreover, as in all previous sections, we will assume that $S$ contains at least one irreducible polynomial. In particular, we may assume that this irreducible polynomial is of Type (I) by Corollary \ref{cor:oddp+irre+nonspecial}. Moreover, we will make separate arguments for the case when $S$ contains exactly one irreducible polynomial and the case when $S$ contains at least two irreducible polynomials.
\subsection{Two irreducibles}
In this subsection, we will assume that $S$ contains at least two irreducible polynomials. In particular, by Corollary \ref{cor:oddp+irre+nonspecial} above, it suffices to assume that \emph{both} of the given irreducible polynomials are of Type (I). From here, we need a few very basic facts. These lemmas are straightforward adaptations of ideas above (see, for instance, the proof of Lemma \ref{lem:oddpbound}), and so we omit their proofs. 
\begin{lem}\label{lem:basicodd} Let $t\geq2$ and $d\geq3$. Then $(t^d-1)^d<t^{d^2}-2t^d$. 
\end{lem}
%\begin{proof}
%Note $(t^d-1)^d<t^{d^2}-2t^d \iff t^{d^2}-(t^d-1)^d>2t^d$.
%Thus we want to show: \[(t^d)^d-(t^d-1)^d=(t^d-(t^d-1))((t^d)^{d-1}+(t^d)^{d-2}(t^d-1)+\dots+(t^d-1)^{d-1})\]\[=\sum_{i=0}^{d-1} (t^d)^{d-i-1}(t^d-1)^{i}>2t^d.\]
%In particular, recall that $t\geq2$ and $d\geq3$, so $t^d\geq8$, giving $(t^d)^{d-1}>(t^d)^2>2t^d$. Furthermore, since $t^d\geq8$, $(t^d)^{d-i-1}(t^d-1)^{i}>0$ for all values of $i$, thus 
%\[\sum_{i=0}^{d-1} (t^d)^{d-i-1}(t^d-1)^{i}=(t^d)^d-(t^d-1)^d>(t^d)^{d-1}>2t^d\]
%\[\Longrightarrow (t^d-1)^d<t^{d^2}-2t^d\] as desired.
%\end{proof}
\begin{lem}\label{lem:basicodd2} For $d\geq2$, the equation $y^d-2=z^d$ has no solutions $y,z\in\mathbb{Z}$ with $|y|\geq2$. 
\end{lem}
%\begin{proof}
%
%\textbf{Case(1):} Suppose $y\in[2,\infty)$.
%Observe that for all ${y}$, $y^{d-n-1}\geq(y-1)^{d-n-1}>0$ and $(y-1)^n\geq(y-2)^n>0$ for $1\leq n\leq d-2$, thus \[\sum_{n=0}^{d-1} y^{d-n-1}(y-1)^n >\sum_{n=0}^{d-1} (y-1)^{d-n-1}(y-2)^n \]
%$$\implies$$
%\[y^d-(y-1)^d>(y-1)^d-(y-2)^d.\]
%Applying repeatedly yields $y^d-(y-1)^d>(y-1)^d-(y-2)^d>(y-2)^d-(y-3)^d>...>2^d-1^d>2$. Thus $y^d-(y-1)^d>2$, giving $(y-1)^d<y^d-2<y^d$. This implies that $y^d-2=z^d$ has no solutions for $y,z\in\mathbb{Z}$ with $|y| \geq 2$, as desired. 
    
%\textbf{Case(2):} Suppose $y\in(-\infty,-2]$ and d is odd.

%Let $y=-k$ for $k\in \mathbb{Z^+}$. Then $y^d-(y-1)^d=(-k)^d-(-k-1)^d=(k+1)^d-k^d$. Since $k\in [2,\infty)$, $k^d-(k-1)^d>(k-1)^d-(k-2)^d$ holds for any $k$, by the argument in Case 1: $(k+1)^d-k^d>2 \implies -k^d>-(k+1)^d+2=(-k-1)^d+2 \implies y^d-2>(y-1)^d$. Giving $(y-1)^d<y^d-2<y^d$, and thus $z$ cannot exist.

%\textbf{Case(3):} Suppose $y\in(-\infty,-2]$ and d is even.    
    %If $d\geq4$, by Case 1, we are done.
    %Consider $d=2$, then let $y=-k$ for $k\in \mathbb{Z^+}$, then we claim %$(-k+1)^2<(-k)^2-2$.
   %Since $k\geq2$, $-2k+1\leq-3<-2 \implies k^2-2k+1<k^2-2 \implies %(-k+1)^2<(-k)^2-2$ which implies $(y+1)^2<y^2-2$, thus giving %$(y+1)^2<y^2-2<y^2$, showing $z$ cannot exist.
%\end{proof}
Next, we classify the integral solutions to some related Diophantine equations.
\begin{prop}\label{prop:oddp} Let $p$ be an odd prime and suppose that $a,b,s,t\in\mathbb{Z}$ are such that 
\begin{equation}\label{eq:new1}
a^p+s^p-s^{p^2}=t^p\;\;\;\;\text{and}\;\;\;\;b^p+t^p-t^{p^2}=s^p.
\end{equation}
Then either $\min\{|s|,|t|\}\leq1$ or $s=t$. 
\end{prop}
\begin{proof} Suppose that $a,b,s,t\in\mathbb{Z}$ satisfy \eqref{eq:new1} and that $s$ and $t$ both have absolute value at least $2$. We consider first the case where $s$ and $t$ have the same parity. In particular, by multiplying both equations in \eqref{eq:new1} by minus one, we can assume that $s$ and $t$ are both positive. In this case, suppose for a contradiction that $t^p>s^p>0$.
Then Lemma \ref{lem:basicodd} implies that 
\begin{equation}\label{eq:new2}
(t^p-1)^p<t^{p^2}-2t^p<t^{p^2}-t^p+s^p=t^{p^2}-(t^p-s^p)<(t^p)^p.
\end{equation}
In particular, $t^{p^2}-t^p+s^p$ is strictly in between consecutive $p$th powers and so cannot be a $p$th power. However, this contradicts the right side of \eqref{eq:new1}. Moreover by symmetry, we see that $s^p>t^p>0$ is also impossible (this time by using the left side of \eqref{eq:new1}). Therefore, we deduce that $s^p=t^p$ and so $s=t$ whenever $s$ and $t$ have the same parity. On the other hand suppose that $s$ and $t$ have mixed parity; without loss say $t$ and $-s$ are both positive. However, if $t>-s$, then the inequalities in \eqref{eq:new2} still hold, and we reach a contradiction. Likewise if $s':=-s>t>0$, then Lemma \ref{lem:basicodd} implies that 
\[((s')^p-1)^p<(s')^{p^2}-2(s')^p<(s')^{p^2}-(s')^p+(-t)^p=(s')^{p^2}-((s')^p+t^p)<((s')^p)^p.\]
In particular, if we set $t'=-t$, then  we see from above that $(s')^{p^2}-(s')^p+(t')^p$ is strictly in between consecutive $p$th powers and so cannot be a $p$th power. However, this contradicts the left side of \eqref{eq:new1}, after we multiply each side by minus one. Therefore, we deduce that $t=-s$ in the case when $s$ and $t$ have mixed parity. On the other hand when $t=-s$, then $a^p=s^{p^2}-2s^p=s^p(s^{p^2-p}-2)$. In particular, since $s\neq0$ it must be the case that $s^{p^2-p}-2=z^p$ for some $z\in\mathbb{Z}$. However, setting $y:=s^{p-1}$ we deduce that $y^p-2=z^p$, which contradicts Lemma \ref{lem:basicodd2}. In particular, when $s$ and $t$ both have absolute value at least $2$, it must be the case that $s=t$. The claim follows. 
\end{proof}
In particular, we deduce the following corollary for sets $S$ containing at least two irreducible polynomials of Type (I). 
\begin{cor}\label{cor:oddp+2irre}
Let $p$ be an odd prime and let $S=\{x^p+c_1,\dots,x^p+c_r\}$ for some distinct $c\in\mathbb{Z}$. Moreover, assume that $S$ contains two irreducible, Type (I) polynomials; call them $\phi_1$ and $\phi_2$. Then either 
\[\{\phi_1^3\circ\phi_2\circ f\,:\, f\in M_S\}\;\;\;\text{or}\;\;\; 
\{\phi_2^3\circ\phi_1\circ f\,:\, f\in M_S\}
\]
is a set of irreducible polynomials. 
\end{cor}
\begin{proof} 
Write $\phi_1(x)=x^p+s^p-s^{p^2}$ and $\phi_2(x)=x^p+t^p-t^{p^2}$ for some $s,t\in\mathbb{Z}$. Note that $s\neq t$ since the maps in $S$ are assumed to be distinct. Moreover, $\phi_1^3$ and $\phi_2^3$ are both irreducible by \cite[Corollary 30]{JonesIMRN}. In particular, if there exist $f,g\in M_S$ such that $\phi_1^3\circ \phi_2\circ f$ and $\phi_2^3\circ \phi_1\circ g$ are reducible polynomials, then Proposition \ref{prop:stability} implies that $\phi_1^3(\phi_2(b))=y^p$ and $\phi_2^3(\phi_1(a))=z^p$ for some $a,b,y,z\in\mathbb{Z}$. But then 
\[a^p+s^p-s^{p^2}=t^p\;\;\;\;\text{and}\;\;\;\;b^p+t^p-t^{p^2}=s^p\]
by Theorem \ref{thm:pthpower}. However, it then follows from Proposition \ref{prop:oddp} that either $s=0,\pm{1}$ or $t=0,\pm{1}$. But then $\phi_1=x^p$ or $\phi_2=x^p$, and we reach a contradiction. Therefore, either $\{\phi_1^3\circ\phi_2\circ f: f\in M_S\}$ or $\{\phi_2^3\circ\phi_1\circ f: f\in M_S\}$ is a set of irreducible polynomials as desired in this case.
\end{proof}
\subsection{One irreducible}
In this subsection, we will assume that $S$ contains an irreducible, Type (I) polynomial $\phi_1=x^p+s^p-s^{p^2}$ and another reducible polynomial $\phi_2$, necessarily of the form $\phi_2(x)=x^p+t^p$ for some $t\in\mathbb{Z}$; see \cite[Corollary 30]{JonesIMRN}. Perhaps interestingly, we need Fermat's famous last theorem to settle this case.   
\begin{thm}[Fermat's Last Theorem \cite{FLT}]\label{FLT} Let $x,y,z\in\mathbb{Z}$ and $n\geq3$. If $x^n+y^n=z^n$, then $xyz=0$. 
\end{thm}
In particular, we have the following consequence for unicritically generated semigroups. 
\begin{cor}\label{oddp+irred+reduc+s notequal t} Let $p$ be an odd prime and let $S=\{x^p+c_1,\dots,x^p+c_r\}$ for some distinct $c_i\in\mathbb{Z}$. Moreover, assume that $S$ contains an irreducible polynomial  $\phi_1(x)=x^p+s^p-s^{p^2}$ of Type (I) and another reducible polynomial $\phi_2(x)=x^p+t^p$ for some $s,t\in\mathbb{Z}$. Then the following statements hold: \vspace{.1cm} 
\begin{enumerate}
\item If $t$ is non-zero and $s\neq t$, then $\{\phi_1^3\circ\phi_2\circ f: f\in M_S\}$ is a set of irreducible polynomials.\vspace{.1cm}  
\item If $t=0$, then there exists $n\gg0$ such that $\{\phi_1^3\circ\phi_2^n\circ f: f\in M_S\}$ is a set of irreducible polynomials.   
\end{enumerate} 
\end{cor}
\begin{proof} Suppose that $S$ contains an irreducible polynomial $\phi_1(x)=x^p+s^p-s^{p^2}$ with a $p$th powered fixed point and another polynomial $\phi_2=x^p+t^p$ that is reducible. Then $\phi_1^3$ is irreducible by \cite[Corollary 30]{JonesIMRN}, and if there exist $f\in M_S$ such that $\phi_1^3\circ \phi_2\circ f$ is reducible, then $\phi_1^3(\phi_2(a))=y^p$ for some $a,y\in\mathbb{Z}$. In particular, Theorem \ref{thm:pthpower} implies that $a^p+t^p=\phi_2(a)=s^p$. Hence, Fermat's Last Theorem \ref{FLT} implies that either $a=0$ or $t=0$; here we use also that $s\neq0$ since $\phi_1$ is irreducible. In particular, if $t$ is non-zero and $t\neq s$ (so $t^p\neq s^p$ also), then we have that $\{\phi_1^3\circ\phi_2\circ f: f\in M_S\}$ is a set of irreducible polynomials as claimed in statement (1). 

Likewise, the case when $t=0$ can be handled similarly: $\phi_1^3\circ \phi_2^n(0)=\phi_1^3(0)$ is not a $p$th power by the proof of \cite[Corollary 30]{JonesIMRN}. Hence, if $\phi_1^3\circ \phi_2^n\circ f$ is reducible, then Proposition \ref{prop:stability} implies that $\phi_1^3(\phi_2^n(a))=y^p$ for some $a,y\in\mathbb{Z}$. However, in this case Theorem \ref{thm:pthpower} implies that $\phi_2^n(a)=a^{p^{n}}=s^p$ and so $a^{p^{n-1}}=s$. But this cannot be true for all $n$ since $s$ is non-zero. Thus, $\{\phi_1^3\circ\phi_2^n\circ f: f\in M_S\}$ is a set of irreducible polynomials for all $n\gg0$ as claimed in statement (2).   
\end{proof}
Hence, combining Corollary \ref{cor:oddp+irre+nonspecial}, Corollary \ref{cor:oddp+2irre}, and Corollary \ref{oddp+irred+reduc+s notequal t}, we answer Question \ref{quest} from the introduction in the affirmative for sets of the form $S=\{x^p+c_1,\dots,x^p+c_r\}$, where $p$ is an odd prime. However, there is a catch: we must assume that is not of the very special form, $S=\{x^p+t^p-t^{p^2},x^p+t^p\}$ for some $t\in\mathbb{Z}$.  

\begin{cor}\label{cor:oddpsummary, outside of special case}
Let $p$ be an odd prime and let $S=\{x^p+c_1,\dots,x^p+c_r\}$ for some distinct $c_i\in\mathbb{Z}$. Moreover, assume that $S$ contains at least one irreducible polynomial and that $S\neq\{x^p+t^p-t^{p^2},x^p+t^p\}$ for all $t\in\mathbb{Z}$. Then $M_S$ contains (an explicit) subset of irreducible polynomials of positive proportion. 
\end{cor}
\begin{proof} Let $\phi_1$ be an irreducible polynomial in $S$. If $\phi_1$ is not of Type (I), then $M_S$ contains a positive proportion of irreducibles by Corollary \ref{cor:irre&non-special}. Therefore, we may assume that $\phi_1(x)=x^p+s^p-s^{p^2}$ for some $s\in\mathbb{Z}$. Likewise, if $S$ contains a second irreducible polynomial $\phi_2$, then we may assume that $\phi_2(x)=x^p+t^p-t^{p^2}$ for some $t\in\mathbb{Z}$. In particular, $M_S$ contains a positive proportion of irreducible polynomials by Corollary \ref{cor:oddp+2irre} in this case. 

Therefore, we may assume that \emph{any} other polynomial in $S$ is reducible. Now take any $\phi_2\in S$ and write $\phi_2(x)=x^p+t^p$ for some $t\in\mathbb{Z}$. Then, if $t=0$ or $t\neq 0$ and $s\neq t$, we have know that $M_S$ contains a positive proportion of irreducibles by Corollary \ref{oddp+irred+reduc+s notequal t}. Hence, we may assume that $x^p+t^p-t^{p^2}$ and $x^p+t^p$ are both elements of $S$. On the other hand, if $\phi_3=x^p+r^p\in S$ is any third polynomial, then we are again done by Corollary \ref{oddp+irred+reduc+s notequal t} (this time applied to $\phi_1$ and $\phi_3$) unless $r=t$. But this implies that $\phi_2=\phi_3$, a contradiction. Therefore, we have succeeded in proving that $M_S$ contains a positive proportion of irreducibles unless $S=\{x^p+t^p-t^{p^2},x^p+t^p\}$ for some $t\in\mathbb{Z}$ as claimed.           
\end{proof} 
\subsection{Remaining special case:} The strategy for attacking the sole remaining case when $S=\{x^p+t^p-t^{p^2},x^p+t^p\}$ for some $t\in\mathbb{Z}$ is as follows. First, we will show that $\phi_1\circ\phi_2\circ\phi_1(a)$ cannot be a $p$th power for all $a\in\mathbb{Z}$; see Lemma \ref{lem:lastcase,refinement} below. In particular, if we can show that the specific polynomial, 
\[\phi_1\circ\phi_2=(x^p+t^p)^p+t^p-t^{p^2},\] 
is irreducible, then it follows that $\{\phi_1\circ\phi_2\circ\phi_1\circ f: f\in M_S\}$ is a set of irreducible polynomials (of positive proportion in $M_S$). However, unfortunately we only succeed in proving that $\phi_1\circ\phi_2$ is irreducible for $p=3$ using some heavy machinery from arithmetic geometry; in fact, we succeed in proving that $\phi_1\circ\phi_2$ is irreducible for all \emph{rational} $t\in\mathbb{Q}$ in this case. Nevertheless, the $p>3$ case remains open (though all data gathered from ``small'' examples supports the claim in general). Moreover, a global technique is necessary to prove this: $\phi_1\circ\phi_2(0)=t^p$ is a $p$th power, and so $\phi_1\circ\phi_2$ is reducible modulo every prime; see Proposition \ref{prop:stability} and Remark \ref{rem:finitefield}. With that said, we continue on with the strategy outlined above.          
\begin{lem}\label{lem:lastcase,refinement}
    Let $\phi_1(x)=x^p+t^p-t^{p^2}$ and $\phi_2(x)=x^p+t^p$ for some odd prime $p$ and $t\in\mathbb{Z}\setminus\{0,\pm{1}\}$. Then, $\phi_1 \circ \phi_2\circ\phi_1(a) \neq y^p$ for all $a,y \in\mathbb{Z}$. 
\end{lem}
\begin{proof}
Assume for the sake of contradiction that $\phi_1 \circ \phi_2\circ\phi_1(a)=y^p$ for some $a,y \in\mathbb{Z}$. Suppose first that $|a|\geq |t|^p+1$. Then 
\[|\phi_1(a)|\geq |a|^p-|t|^{p^2}-|t|^p\geq (|t|^p+1)^p-|t|^{p^2}-|t|^p\geq p|t|^{p(p-1)}.\]
In particular, we similarly deduce that 
\[|\phi_2(\phi_1(a))|\geq |\phi_1(a)|^p-|t|^p>p^p|t|^{p^2(p-1)}-|t|^p>(p^p-1)|t|^{p^2(p-1)}.\]
On the other hand, Lemma \ref{lem:oddpbound}, the fact that $\phi_1 \circ \phi_2\circ\phi_1(a)=y^p$, and the bound above together imply that 
\[
\scalemath{.9}{
(p^p-1)|t|^{p^2(p-1)}\leq|\phi_2(\phi_1(a))|\leq \Big(\frac{|t^{p^2}-t^p|}{p}\Big)^{1/(p-1)}+1\leq
%\Big(\frac{|t|^{p^2}+|t|^p}{p}\Big)^{1/(p-1)}\leq 
\Big(\frac{2|t|^{p^2}}{p}\Big)^{1/(p-1)}+1\leq |t|^{p^2/(p-1)}+1.}\]
However, this is clearly false unless $t=0$, and we reach a contradiction. Likewise, when $|a|<|t|^p-1$ we have that 
\begin{align*} 
|\phi_1(a)|=|-\phi_1(a)|&=|t^{p^2}-a^p-t^p|\geq|t|^{p^2}-{|a|}^{p}-|t|^p\\[3pt] 
&>\Big(|t|^{p^2}-(|t|^p-1)^p\Big)-|t|^p>p(|t|^p-1)^{p-1}-|t|^p\\[3pt] 
&>(p-1)(|t|^p-1)^{p-1}.
\end{align*}
In particular, we deduce that \vspace{.1cm} 
\[
\scalemath{.93}{
|\phi_2(\phi_1(a))|\geq |\phi_1(a)|^p-|t|^p>(p-1)^p(|t|^p-1)^{p(p-1)}-|t|^p>((p-1)^p-1)(|t|^p-1)^{p(p-1)}.}\vspace{.1cm}\]
However, as before the lower bound above clearly exceeds the upper bound
\[|\phi_2(\phi_1(a))|\leq \Big(\frac{|t^{p^2}-t^p|}{p}\Big)^{1/(p-1)}+1\leq
%\Big(\frac{|t|^{p^2}+|t|^p}{p}\Big)^{1/(p-1)}\leq 
\Big(\frac{2|t|^{p^2}}{p}\Big)^{1/(p-1)}+1\leq |t|^{p^2/(p-1)}+1 \vspace{.1cm} \]
coming from Lemma \ref{lem:oddpbound} and the triangle inequality. Therefore, we must have that $|a|=|t|^p$ so that $a=\pm{t^p}$. However, $|\phi_1(t^p)|=|t|^p$ and $|-\phi_1(-t^p)|\geq2|t|^{p^2}-|t|^p$. Therefore, if $a=\pm{t^p}$, then we see that  
\[ |\phi_1(a)|\geq\min\{|t|^p,2|t|^{p^2}-|t|^p\}=|t|^p.\]
But then the bound from Lemma \ref{lem:oddpbound} implies that 
\[|t|^{p^2}-|t|^p\leq|\phi_1(a)|^p-|t|^p\leq |\phi_2(\phi_1(a))|\leq |t|^{p^2/(p-1)}+1.\]
However, it then follows that 
\[|t|^{p^2}\leq  |t|^{p^2/(p-1)}+|t|^p\leq 2 |t|^{p^2/(p-1)}\leq |t|^{\frac{p^2}{p-1}+1}+1\]
for all $|t|\geq2$. But $p^2>p^2/(p-1)+1$ for all $p\geq3$, and we reach a contradiction. The result follows.      
\end{proof}
Next, we now show that $\phi_1\circ\phi_2$ is irreducible for all suitable $t\in\mathbb{Q}$ when $p=3$. 
\begin{prop}\label{prop:cubic,specialcase} Let $\phi_1(x)=x^3+t^3-t^9$ and $\phi_2(x)=x^3+t^3$ for some $t\in\mathbb{Q}\setminus\{0,\pm{1}\}$. Then $\phi_1\circ\phi_2$ is irreducible in $\mathbb{Q}[x]$.  
\end{prop}
\begin{proof} Let $\zeta$ be a primitive $3$rd root of unity, let $L=\mathbb{Q}(\zeta)$, and let $t\in\mathbb{Q}\setminus\{0,\pm{1}\}$. We first show that $\phi_1=x^3+t^3-t^{9}$ is irreducible in $L[x]$. To do this, we note that $\phi_1(0)=t^3-t^9$ is not a cube $\mathbb{Q}$. Indeed, if $z^3=t^3-t^9=t^3(1-t^6)$ for some $z\in\mathbb{Q}$, then the fact that $t$ is non-zero implies that $y^3=1-t^6$ for some $y\in\mathbb{Q}$. In particular, we obtain a rational point $(X,Y)=(-y,t^3)$ on the elliptic curve $E: Y^2=X^3+1$. However, we compute with {\tt{Magma}} \cite{Magma} that $E(\mathbb{Q})\cong\mathbb{Z}/6\mathbb{Z}$ and $E(\mathbb{Q})=\{\infty,(0,\pm{1}), (-1,0),(2,\pm{3})\}$. In particular, $t^3\in\{0,\pm{1},\pm{3}\}$, which is impossible since $t\in\mathbb{Q}\setminus\{0,\pm{1}\}$. Therefore, $\phi_1$ has no roots in $\mathbb{Q}$ and so is irreducible in $\mathbb{Q}[x]$. However, this implies that $\phi_1$ has no roots in $L$ either: if $\alpha\in\overline{\mathbb{Q}}$ is a root of $\phi_1$, then $[\mathbb{Q}(\alpha):\mathbb{Q}]=3$ so that $L\cap \mathbb{Q}(\alpha)=\mathbb{Q}$ follows for the fact that $[L:\mathbb{Q}]=2$. Therefore, $\phi_1$ is irreducible over $L$ as claimed. 

Now, assume for a contradiction that $\phi_1\circ\phi_2$ is reducible in $\mathbb{Q}[x]$. Then $\phi_1\circ\phi_2$ is reducible in $L[x]$ and so  \cite[Theorem 10]{JonesIMRN} implies that 
\begin{equation}\label{eq:cyclotomicfactorization}
\phi_1\circ\phi_2(x)=h(x)\cdot h(\zeta x)\cdot h(\zeta^{2}x)
\end{equation} 
for some monic, irreducible polynomial $h\in L[x]$.
% so the other factors h(zeta^k) are irreducible in L[x] too
In particular, we see that $\deg(h)=3$. Hence, it follows that any irreducible factor $g\in\mathbb{Q}[x]$ has degree divisible by $3$. This leaves two possibilities: either $\phi_1\circ\phi_2$ factors in $\mathbb{Q}[x]$ as the product of an irreducible cubic and an irreducible sextic or the product of three irreducible cubics (i.e, $9=6+3$ and $9=3+3+3$ are the only partitions of $9$ into parts divisible by $3$). In particular, we may assume that $\phi_1\circ\phi_2$ has a monic irreducible factor $g\in\mathbb{Q}[x]$ of degree $3$. Note that such a $g$ must also be irreducible in $L[x]$. Hence, by unique factorization in $L[x]$, we have that $g(x)=h(\zeta^k x)$ for some $k$. In particular, by replacing $h(x)$ with $h(\zeta^k x)$ if necessary, we may assume that $h\in\mathbb{Q}[x]$. Note that \eqref{eq:cyclotomicfactorization} implies that $t^3=\phi_1\circ\phi_2(0)=h(0)^3$ and so $h(0)=t$. Hence, $h(x)=x^3+ax^2+bx+t$ for some $a,b\in\mathbb{Q}$ and 
\begin{equation*}
x^9+(3t^3)x^6+(3t^6)x^3+t^3=x^9 + (a^3 - 3ab + 3t)x^6 + (-3abt + b^3 + 3t^2)x^3 + t^3    
\end{equation*}
follows from \eqref{eq:cyclotomicfactorization}. Moreover, by equating coefficients, eliminating $b$
% by solving for b in the x^6 term and subing it into the x^3 term
and making the substitution $a^3=s$, we see that
\begin{equation}\label{eq:cubic,genus2}
B:\, s^3 - 9s^2t^3 - 18s^2t - 54st^6 + 27st^4 + 27st^2 - 27t^9 + 81t^7 - 81t^5 + 27t^3=0
\end{equation}
for some $s,t\in\mathbb{Q}$. Moreover, we compute with {\tt{Magma}} \cite{Magma} that the genus 2 curve above is birational over $\mathbb{Q}$ to the hyperelliptic curve:
\begin{equation}\label{eq:cubic,hyperelliptic}
C:\; Y^2=-X^6 + 3X^5 + 3X^4 - 11X^3 + 3X^2 + 3X - 1.
\end{equation}
Furthermore, we find six rational points $(X,Y)=(-1,\pm{3}), (2,\pm{3}, (1/2,\pm{3/8})$ which correspond to the points on \eqref{eq:cubic,genus2} with $t$-coordinates in $\{0,\pm{1}\}$ and a point at infinity. Therefore, it suffices to show that $C(\mathbb{Q})$ consists exactly of the six known rational points to prove the proposition. On the one hand, the Jacobian of $C$ has rank two (equal to the genus of $C$), and so the usual Chabauty-Coleman method does not directly apply. However, the Galois group of the polynomial $-X^6 + 3X^5 + 3X^4 - 11X^3 + 3X^2 + 3X - 1$ is quite small, isomorphic to $\mathbb{Z}/3\mathbb{Z}$, and so we can hope that a blending of other techniques (e.g., unramified covering collections and the elliptic Chabuaty method \cite{Bruin,stoll2011rational}) can yield the desired conclusion. This is indeed the case, which we now discuss in some detail.

We begin by factoring,
\begin{align*}
-X^6 + 3X^5 + 3X^4 - 11X^3 + 3X^2 + 3X - 1 =&(X^4 - \theta X^3 - 3X^2 + (3\theta + 1)X - \theta)\\[2pt]
&\;\;\;\;\;\;\cdot(-X^2 + (-\theta + 3)X - \theta^2 + 3\theta),
\end{align*}
where $\theta$ is a root of the irreducible cubic $X^3 - 3X^2 + 1$. Let $K=\mathbb{Q}(\theta)$, a number field with class number equal to one and unit group of rank two. In particular, it follows from \cite[\S2.3]{stoll2011rational} that each rational point $(X,Y)\in C(\mathbb{Q})$ lifts to a rational point $(X,U,V)\in\mathcal{D}_d(\mathbb{Q})$, where 
\[
\scalemath{.95}{
\mathcal{D}_d\,: dU^2=X^4 - \theta X^3 - 3X^2 + (3\theta + 1)X - \theta\;\;\text{and}\;\; dV^2=-X^2 + (-\theta + 3)X - \theta^2 + 3\theta}
\]
for some $d\in K^*/(K^*)^2$ supported on the primes dividing the resultant of of the polynomials $F(X)=X^4 - \theta X^3 - 3X^2 + (3\theta + 1)X - \theta$ and $G(X)=-X^2 + (-\theta + 3)X - \theta^2 + 3\theta$; formally, the the maps $\mathcal{D}_d(\mathbb{Q})\rightarrow C(\mathbb{Q})$ given by $(X,U,V)\rightarrow (X,dUV)$ are unramified $\mathbb{Z}/2\mathbb{Z}$-coverings and $\mathcal{D}_d$ has points everywhere locally if and only if $d$ is supported on the primes dividing $\text{Res}(F,G)=27\theta^2 - 54\theta - 81$. In particular, there are a possible $16$ choices (up to square class) for $d$: 
\begin{equation*}\label{eq:twists}
d\in T:=\Big\{(-1)^{e_1}\cdot(-\theta^2 + 2\theta + 1)^{e_2}\cdot \theta^{e_3}\cdot (-1-3\theta)^{e_4}\;:\; e_1,e_2,e_3,e_4\in\{0,1\}\Big\}.
\end{equation*}
On the other hand, the equations defining the left side of $\mathcal{D}_d$ are elliptic curves, 
\[E_d:\, dU^2=X^4 - \theta X^3 - 3X^2 + (3\theta + 1)X - \theta\]
defined over $K$. Moreover, when $\rank(E_d(K))<[K:\mathbb{Q}]$, then one can often \cite{Bruin} use the formal group on $E_d$ at some auxiliary prime together with bounds on the number of roots of certain power series to describe the \emph{finite sets}, 
\[F_d:=\{X\in\mathbb{Q}\,:\, (X,U)\in E_d(K)\;\text{for some}\; U\in K\}.\]
In fact, this procedure has (largely) been implemented in {\tt{Magma}} \cite{Magma} (except for the fact that we must first compute a Weierstrass equation for each $E_d$, but this may be done with {\tt{Magma}} also). We carry out this procedure and find that $\text{rank}(E_d(K))<2$ and 
\[\bigcup_{d\in T}F_d=\{0, \pm{1}, 1/2, 2, \infty\}.\]
However, among these only $X=-1,1/2,2$ yield rational points on $C$. Therefore, we have found all of the rational points on $C$. In particular, we deduce that if $\phi_1\circ\phi_2$ is reducible in $\mathbb{Q}[x]$, then $t=0,\pm{1}$ as claimed; see \href{https://sites.google.com/a/alumni.brown.edu/whindes/research/code}{Mathworks 2023} for {\tt{Magma}} \cite{Magma} code. 
\end{proof}
In particular, we have the following immediate consequence: 
\begin{cor}\label{cor:p=3}
Let $S=\{x^3+t^3-t^9, x^3+t^3\}$ for some $t\in\mathbb{Z}\setminus\{0,\pm{1}\}$. Then $\{\phi_1\circ\phi_2\circ\phi_1\circ f: f\in M_S\}$ is a set of irreducible polynomials. 
\end{cor}
\begin{proof}
Combine Proposition \ref{prop:cubic,specialcase}, Lemma \ref{lem:lastcase,refinement}, and Proposition \ref{prop:stability}.     
\end{proof}
\subsection{Summary and Loose Ends}
We end this section by summarizing our findings in the case of odd primes in a single theorem. This may be viewed as an explicit version of Theorem \ref{thm:intro-oddprime} from the introduction; compare to our main result Theorem \ref{thm:main} in the quadratic case. Lastly, we give a short proof of our local-global result from the introduction; see Theorem \ref{thm:local-global}.    
\begin{thm}\label{thm:oddprime,explicit} Let $S=\{x^p+c_1,\dots,x^p+c_r\}$ for some odd prime $p$ and some distinct $c_1,\dots,c_r\in\mathbb{Z}$. Moreover, assume that $S$ contains an irreducible polynomial. Then the following statements hold:\vspace{.15cm}
\begin{enumerate}
\item If $r=1$, then $M_S$ is a set of irreducible polynomials. \vspace{.2cm}
\item If $S$ contains an irreducible polynomial $\phi_1$ that does not have a $p$th powered fixed point, then 
\[\{\phi_1^3\circ f: f\in M_S\}\]
is a set of irreducible polynomials. \vspace{.2cm}
\item If $S$ contains two irreducible polynomials $\phi_1$ and $\phi_2$ that both have $p$th powered fixed points, then either 
\[\{\phi_1^3\circ\phi_2\circ f: f\in M_S\}\qquad \text{or}\qquad \{\phi_2^3\circ\phi_1\circ f: f\in M_S\}\]
is a set of irreducible polynomials.\vspace{.2cm} 
\item If $S$ contains an irreducible polynomial $\phi_1(x)=x^p+s^p-s^{p^2}$ with a $p$th powered fixed point and another reducible polynomial $\phi_2(x)=x^p+t^p$ for some non-zero $t$ with $t\neq s$, then 
\[\{\phi_1^3\circ \phi_2\circ f: f\in M_S\}\]
is a set of irreducible polynomials.\vspace{.2cm}
\item If $S$ contains an irreducible polynomial $\phi_1(x)=x^p+s^p-s^{p^2}$ with a $p$th powered fixed point and another polynomial $\phi_2(x)=x^p$, then there exist $n\gg0$ such that 
\[\{\phi_1^3\circ \phi_2^n\circ f: f\in M_S\}\]
is a set of irreducible polynomials.
 \item If $p=3$ and $S$ contains an irreducible polynomial  $\phi_1(x)=x^p+s^p-s^{p^2}$ with a $p$th powered fixed point and another reducible polynomial of the form $\phi_2(x)=x^p+s^p$, then \[\{\phi_1\circ \phi_2\circ\phi_1\circ f: f\in M_S\}\]
is a set of irreducible polynomials.
\vspace{.1cm}
\end{enumerate}
In particular, if $p=3$ then $M_S$ contains a positive proportion of irreducible polynomials if and only if it contains at least one irreducible polynomial. More generally, if $p>3$ and $S\neq\{x^p+t^p-t^{p^2}, x^p+t^p\}$ for all $t\in\mathbb{Z}$, then $M_S$ contains a positive proportion of irreducible polynomials if and only if it contains at least one irreducible polynomial 
\end{thm}
\begin{proof} Statements (1)-(6) follow from \cite[Corollary 30]{JonesIMRN}, Corollary \ref{cor:oddp+irre+nonspecial}, Corollary \ref{cor:oddp+2irre}, Corollary \ref{oddp+irred+reduc+s notequal t}, and Corollary \ref{cor:p=3} respectively.   
\end{proof}
Now for the proof of the local-global result from the introduction. As a reminder, we let $p=2$ or $p=3$, let $t\in\mathbb{Z}\setminus\{0,\pm{1}\}$, and let  
\begin{equation}\label{local-global} 
S=\big\{x^p+t^p-t^{p^2},\;x^p+(-1)^{p-1}t^p,\;x^p+c_1,\;\dots,\;x^p+c_r\big\}
\end{equation} 
for some arbitrary $c_1,\dots,c_r\in\mathbb{Z}$. Then we will show that the  polynomials 
\[I_2:=\{\phi_1^2\circ\phi_2\circ\phi_1\circ f: f\in M_S\}\;\;\;\text{and}\;\;\;  
I_3:=\{\phi_1\circ\phi_2\circ\phi_1\circ f:f\in M_S\}\] 
break the local-global principle for irreducibility when $p=2$ and $p=3$ respectively.  
\begin{proof}[(Proof of Proposition \ref{thm:local-global})] We already know that $I_2$ and $I_3$ are sets of irreducible polynomials (in $\mathbb{Q}[x]$) by Theorem \ref{thm:family1+s=t} and Corollary \ref{cor:p=3} respectively. As for local information, we consider the $p=2$ and $p=3$ cases separately. 

When $p=2$, we have that $\phi_1^2\circ \phi_2(0)=t^2$ and thus $\phi_1^2\circ \phi_2(0)$ is a square modulo $q$ for every prime $q$. In particular, $\phi_1^2\circ \phi_2$ is reducible in $\mathbb{F}_q[x]$ for every odd prime $q$ by Proposition \ref{prop:stability} and Remark \ref{rem:finitefield}; hence, every polynomial of the form $\phi_1^2\circ \phi_2\circ\phi_1\circ f$ is likewise reducible modulo $q$. On the other hand, $\phi_1(x)\equiv (x+t-t^2)^2 \mod{2}$, and so every polynomial of the form $\phi_1^2\circ\phi_2\circ\phi_1\circ f$ is reducible in $\mathbb{F}_2[x]$ as well. Therefore, $I_2$ is a set of polynomials that break the local-global principle for irreducibility as claimed. 

Likewise, when $p=3$ we have that $\phi_1\circ \phi_2(0)=t^3$ and so $\phi_1\circ \phi_2(0)$ is a cube modulo $q$ for every prime $q$. In particular, $\phi_1\circ \phi_2$ is reducible in $\mathbb{F}_q[x]$ for every prime $q\neq3$ by Proposition \ref{prop:stability} and Remark \ref{rem:finitefield}. Therefore, every polynomial of the form $\phi_1\circ \phi_2\circ\phi_1\circ f$ is likewise reducible modulo $q$. On the other hand, $\phi_1(x)\equiv (x+t-t^3)^3 \mod{3}$, and so every polynomial of the form $\phi_1\circ\phi_2\circ\phi_1\circ f$ is reducible in $\mathbb{F}_3[x]$ as well. Hence, $I_3$ is a set of polynomials that break the local-global principle for irreducibility as claimed. 
\end{proof}
Finally, we end this paper with a conjecture that the special polynomials discussed in subsection 4.3 are irreducible for all odd primes $p$ (not just $p=3$).  
\begin{conjecture}
Let $p$ be an odd prime and let $t\in\mathbb{Z}\setminus\{0,\pm{1}\}$. Then the polynomial $(x^p+t^p)^p+t^p-t^{p^2}$ is irreducible in $\mathbb{Q}[x]$. 
\end{conjecture}
%\begin{rem} Let $f_{p,t}(x)=(x^p+t^p)^p+t^p-t^{p^2}$ and let $v_p(\cdot)$ %denote the $p$-adic valuation. Then it is straightforward to check that if %$v_p(t)=0$ and $v_p(t^{p^2-p}-1)\not\equiv 0\mod p$, then the shifted %polynomial $f_{p,t}(x-t)$ satisfies Dumas's irreduciblity criterion. In particular, it follows that $f_{p,t}(x)$ is irreducible for many values of $p$ and $t$.       
%\end{rem}
In particular, provided that the conjecture above is true, it follows that for all 
\[S=\{x^p+c_1,\dots,x^p+c_r\}\] 
and $c_i\in\mathbb{Z}$, the semigroup $M_S$ contains a positive (and explicit) proportion of irreducible polynomials if and only if it contains at least one irreducible polynomial. In particular, one may remove the only remaining case corresponding to $S=\{x^p+t^p-t^{p^2},x^p+t^p\}$. 
\begin{rem} Let $S$ be as in \eqref{local-global} for some odd prime $p$. Then the conjecture above implies that the polynomials in $\{\phi_1\circ\phi_2\circ\phi_1\circ f\,:\, f\in M_S\}$ break the local-global principle for irreducibility (as is known for $p=3$). 
\end{rem}
\section{Results over the rationals}
In this section, we adapt ideas from previous sections (and from \cite{Mathworks}) to produce irreducible polynomials in semigroups with \emph{rational} coefficients. However, our results are more restrictive in this case. In particular, we need the following definition \cite[\S5]{jones2013galois}.  
\begin{defn} Let $\phi(x)\in\mathbb{Q}[x]$. Then we say that $\phi$ is stable (over $\mathbb{Q}$) if $\phi^n$ is irreducible in $\mathbb{Q}[x]$ for all $n\geq1$.     
\end{defn}
\begin{rem} In particular, if $\phi(x)=x^p+c$ for some $c\in\mathbb{Z}$ and $p$ is prime, then $\phi$ is stable if and only if $\phi$ is irreducible; see \cite[Corollary 30]{JonesIMRN} and \cite[Corollary 1.3]{Stoll}. However, this is not true in general for rational values of $c$. For some results on stable polynomials (and other similar properties, e.g., eventual stability), see \cite{Sadek,MR4103000,jones2012iterative,jones2017eventually}.       
\end{rem}
Likewise, we need some assumptions about the possible periodic points for the maps in our generating set. Recall that a point $a\in\mathbb{Q}$ is periodic for $\phi\in\mathbb{Q}[x]$ if $\phi^n(a)=a$ for some $n\geq1$. Moreover, the minimum such $n$ is called the \emph{period length} of $a$. In particular, when $S=\{x^2+c_1,\dots,x^2+c_r\}$ we will assume the following famous conjecture in arithmetic dynamics; see \cite{flynn1997cycles,hutz2013poonen,Morton,Poonen,stoll2008rational} for some partial results.     
\begin{conjecture}[Poonen \cite{Poonen}] Let $\phi(x)=x^2+c$ for some $c\in\mathbb{Q}$. If $a\in\mathbb{Q}$ is a periodic point for $\phi$, then the period length of $a$ is at most $3$.    
\end{conjecture} 
Likewise, if $\phi^m(a)$ periodic for $\phi$ for some $m\geq0$, then we say that $a$ is \emph{preperiodic} for $\phi$. Moreover, we let $\Prep(\phi)$ denote the set of all preperiodic points for $\phi$. 

To prove our result over the rational numbers, we need a few facts. First, we prove that if $\phi(x)=x^p+c$ for $c\in\mathbb{Q}$ is irreducible and there exists a large iterate $n$ and rational numbers $a,y\in\mathbb{Q}$ such that $\phi^n(a)=y^p$, then $y^p$ must be a periodic point for $\phi$; compare this to more explicit results in Theorem \ref{thm:4th=square} and Theorem \ref{thm:pthpower} when $c$ is integral. Roughly speaking, in this case we use the finiteness of rational points on certain curves to get height bounds on iterates, instead of divisibility properties to get absolute value bounds on iterates.     
\begin{lem}\label{lem:n} Let $\phi(x)=x^p+c$ for some $c\in\mathbb{Q}$ and some $p\geq2$. Moreover, assume that $\phi$ is irreducible in $\mathbb{Q}[x]$. Then there is an $n=n_\phi\gg0$ such that if there exist $a,y\in\mathbb{Q}$ with $\phi^{n}(a)=y^p$, then $y^p$ is a periodic point of $\phi$.        
\end{lem}
\begin{proof} For $p\geq2$, consider the curves with affine models given by 
\begin{equation*}
C_p = 
\left\{
    \begin{array}{lr}
       \phi(x)=y^p, & \text{if } p \geq 4\\[2pt]
        \phi^2(x)=y^3, & \text{if } p=3\\[2pt] 
        \phi^3(x)=y^2, & \text{if } p=2\\
    \end{array}
\right\} 
\end{equation*}
Then for $p\geq4$, the genus of $C_p$ is $(p-1)(p-2)/2>1$ by the Hurwitz genus formula. Likewise, the genus of $C_3$ is $7$, and the genus of $C_2$ is $3$. Hence, $C$ has genus at least 2 for all $p\geq2$, and so $C_p(\mathbb{Q})$ is finite by Faltings Theorem. In particular, there is a bound $B=B_p$ on the height of the $x$-coordinates of the elements of $C_p(\mathbb{Q})$. Now let $\hat{h}_\phi:\mathbb{P}^1\rightarrow\mathbb{R}_{\geq0}$ denote the canonical height function attached to $\phi$ given by the limit $\hat{h}_\phi(a)=\lim_{n\rightarrow\infty} h(\phi^n(a))/\deg(\phi)^n$, where $h(\cdot)$ is the usual logarithmic Weil height; see \cite[\S3.1]{SilvermanBook}. Then $\hat{h}_\phi$ has the following properties: 
\begin{equation}\label{eq:canonicalht} 
|h(a)-\hat{h}_\phi(a)|\leq C_\phi\;\;\;\text{and}\;\;\; \hat{h}_\phi(\phi^n(a))=\deg(\phi^n)\cdot \hat{h}_\phi(a)
\end{equation} 
for some non-negative constant $C_\phi$ and all $a\in\mathbb{P}^1(\overline{\mathbb{Q}})$ and $n\geq1$; see \cite[\S3.4]{SilvermanBook}. In particular, it follows from Northcott's theorem \cite[Theorem 3.12]{SilvermanBook} that $\hat{h}_\phi(a)=0$ if and only $a\in\Prep(\phi)$ and the quantity, 
\begin{equation}\label{hatmin}
\hat{h}_{\phi,\mathbb{Q}}^{\min}:=\inf\{\hat{h}_\phi(a)\,: a\in\mathbb{Q}\;\text{and}\; a\not\in\Prep(\phi)\},
\end{equation} 
is strictly positive; see, \cite{Ingram} for explicit lower bounds. Now suppose that $\phi^n(a)=y^p$ for some $a,y\in\mathbb{Q}$ and that $n>\log_p((B+C_\phi)/\hat{h}_{\phi,\mathbb{Q}}^{\min})+2$. We will show that $a\in\Prep(\phi)$. To see this, note that $(\phi^{n-i}(a),y)\in C_p(\mathbb{Q})$ for some $0\leq i\leq 2$ and so 
\[p^{n-2}\cdot\hat{h}_{\phi,\mathbb{Q}}^{\min}-C_\phi \leq p^{n-i}\cdot\hat{h}_\phi(a)-C_\phi=\hat{h}_\phi(\phi^{n-i}(a))-C_\phi\leq h(\phi^{n-i}(a))\leq B\]
follows from \eqref{eq:canonicalht} and \eqref{hatmin}. However, this contradicts the assumed lower bound on $n$. Therefore, $a\in\Prep(\phi)$ as claimed. However, there are only a finitely many such $a\in\mathbb{Q}$ by Northcotts theoreom \cite[Theorem 3.12]{SilvermanBook} (and the properties of $\hat{h}_\phi$ discussed above). Therefore, we may assume also (by enlarging our lower bound if necessary) that $n$ is larger than the tail lengths of each of the finitely many preperiodic portraits in $\mathbb{Q}$. That is, there is an $n_\phi\gg0$ so that if $\phi^n(a)=y^p$ for some $a,y\in\mathbb{Q}$, then $a\in\Prep(\phi)$ and $\phi^n(a)=y^p$ is a \emph{periodic point} as claimed.        
\end{proof} 

\begin{lem}\label{lem:onlyfixedpts} Let $\phi(x)=x^p+c$ for some $c\in\mathbb{Q}$ and some odd $p\geq3$. Then any rational periodic point of $\phi$ is a fixed point.  
\end{lem}
\begin{proof} This follows from work in \cite{MR3447581}; see also \cite[Theorem 3]{Hutz}.     
\end{proof}

\begin{lem}\label{lem:nosquaresin3cycle} Let $\phi(x)=x^2+c$ for some $c\in\mathbb{Q}$ and assume that $y\in\mathbb{Q}$ is a point of exact period length equal to three. Then $y$ is not a square in $\mathbb{Q}$.  
\end{lem}
\begin{proof} Let $\phi$ and $y$ be as above. Then it follows from \cite[Theorem 3]{Morton} that: 
\[
\scalemath{.89}{
c=-\frac{t^6+2t^5+4t^4+8t^3+9t^2+4t+1}{4t^2(t+1)^2}\;\text{and}\;\;y=\frac{t^3+2t^2+t+1}{2t(t+1)},\frac{t^3-t-1}{2t+1},\;\text{or}\;\frac{t^3+2t^2+3t+1}{2t(t+1)}} 
\]
for some $t\in\mathbb{Q}\setminus\{0,-1\}$; see also \cite[Theorem 1]{Poonen}. In particular, if $y$ is a square then we obtain a rational point on one of the hyperelliptic curves:
\begin{align*} 
B_1&: Y^2=2t(t+1)(t^3+2t^2+t+1), \\[3pt] 
B_2&: Y^2=2t(t+1)(t^3-t-1), \\[3pt]
B_3&: Y^2=2t(t+1)(t^3+2t^2+3t+1).
\end{align*}
In particular, it suffices to show that the only affine points in $B_i(\mathbb{Q})$ correspond to $t=0,-1$. However, these curves each have genus two and rank zero Jacobians, and so the Chabuaty method (which is implemented fully in {\tt{Magam} \cite{Magma}} in this special case) can be used to prove that $B_i(\mathbb{Q})=\{\infty, (0,0), (-1,0)\}$ for all $i$ as desired; see   \href{https://sites.google.com/a/alumni.brown.edu/whindes/research/code}{Mathworks 2023} for the code verifying this calculation.    
\end{proof} 
We now have the tools in place to prove Theoreom \ref{thm:rational} on semigroups generated by sets of the form $S=\{x^p+c_1,\dots,x^p+c_r\}$ for $c_1,\dots,c_r\in\mathbb{Q}$ and $p$ a prime. 
\begin{proof}[(Proof of Theorem \ref{thm:rational})] Let $S$ be as above and assume that some $\phi\in S$ is stable. Moreover, let $n=n_\phi$ be as in Lemma \ref{lem:n} Then, if some polynomial of the form $\phi^n\circ f$ is reducible, Proposition \ref{prop:stability} implies that $\phi^n(a)=y^p$ for some $a,y\in\mathbb{Q}$. But then Lemma \ref{lem:n} implies that $y^p$ is a rational periodic point of $\phi$. In particular, if $p$ is odd, then Lemma \ref{lem:onlyfixedpts} implies that $y^p$ is a fixed point. Hence, $\phi$ is of Type (I) as desired. Likewise if $p=2$,  then Poonen's conjecture implies that $y^2$ must have period length at most $3$. However, Lemma \ref{lem:nosquaresin3cycle} implies that the period length $3$ case cannot occur. Therefore, $\phi$ is of Type (I) or of Type (II) as desired.           
\end{proof}
\bibliographystyle{plain}
\bibliography{Irreducibles}

\begin{thebibliography}{10}

\bibitem{Magma}
Wieb Bosma, John Cannon, and Catherine Playoust.
\newblock The magma algebra system i: The user language.
\newblock {\em Journal of Symbolic Computation}, 24(3-4):235--265, 1997.

\bibitem{brookfield2007factoring}
Gary Brookfield.
\newblock Factoring quartic polynomials: a lost art.
\newblock {\em Mathematics Magazine}, 80(1):67--70, 2007.

\bibitem{Bruin}
Nils Bruin.
\newblock Chabauty methods using elliptic curves.
\newblock {\em J. Reine Angew. Math.}, 562:27--49, 2003.

\bibitem{Octics}
Malcolm Hoong~Wai Chen, Angelina Yan~Mui Chin, and Ta~Sheng Tan.
\newblock Galois groups of certain even octic polynomials.
\newblock {\em Journal of Algebra and Its Applications}, page 2350263, 2022.

\bibitem{Sadek}
Mohamed~O. Darwish and Mohammad Sadek.
\newblock Eventual stability of pure polynomials over the rational field, 2023.
\newblock submitted.

\bibitem{MR4103000}
David DeMark, Wade Hindes, Rafe Jones, Moses Misplon, Michael Stoll, and
  Michael Stoneman.
\newblock Eventually stable quadratic polynomials over {$\mathbb Q$}.
\newblock {\em New York J. Math.}, 26:526--561, 2020.

\bibitem{flynn1997cycles}
EV~Flynn, Bjorn Poonen, and Edward~F Schaefer.
\newblock Cycles of quadratic polynomials and rational points on a genus-$2 $
  curve.
\newblock {\em Duke Math. J.}, 90(1):435--463, 1997.

\bibitem{JonesIMRN}
Spencer Hamblen, Rafe Jones, and Kalyani Madhu.
\newblock The density of primes in orbits of zd+ c.
\newblock {\em International Mathematics Research Notices}, 2015(7):1924--1958,
  2015.

\bibitem{DiscCont}
Wade Hindes.
\newblock Orbit counting in polarized dynamical systems.
\newblock {\em Discrete Contin. Dyn. Syst.}, 42(1):189--210, 2022.

\bibitem{Mathworks}
Wade Hindes, Reiyah Jacobs, and Peter Ye.
\newblock Irreducible polynomials in quadratic semigroups.
\newblock {\em J. Number Theory}, 248:208--241, 2023.

\bibitem{Hutz}
Benjamin Hutz.
\newblock Determination of all rational preperiodic points for morphisms of
  {PN}.
\newblock {\em Math. Comp.}, 84(291):289--308, 2015.

\bibitem{hutz2013poonen}
Benjamin Hutz and Patrick Ingram.
\newblock On poonen's conjecture concerning rational preperiodic points of
  quadratic maps.
\newblock {\em The Rocky Mountain Journal of Mathematics}, pages 193--204,
  2013.

\bibitem{Ingram}
Patrick Ingram.
\newblock Lower bounds on the canonical height associated to the morphism
  {$\phi(z)=z^d+c$}.
\newblock {\em Monatsh. Math.}, 157(1):69--89, 2009.

\bibitem{jones2012iterative}
Rafe Jones.
\newblock An iterative construction of irreducible polynomials reducible modulo
  every prime.
\newblock {\em Journal of Algebra}, 369:114--128, 2012.

\bibitem{jones2013galois}
Rafe Jones.
\newblock Galois representations from pre-image trees: an arboreal survey.
\newblock {\em Publ. Math. Besan\c{c}on Alg\`ebre Th\'{e}orie Nr.}, pages
  107--136, 2013.

\bibitem{jones2017eventually}
Rafe Jones and Alon Levy.
\newblock Eventually stable rational functions.
\newblock {\em International Journal of Number Theory}, 13(09):2299--2318,
  2017.

\bibitem{Morton}
Patrick Morton.
\newblock Arithmetic properties of periodic points of quadratic maps.
\newblock {\em Acta Arith.}, 62(4):343--372, 1992.

\bibitem{MR3447581}
W.~Narkiewicz.
\newblock On a class of monic binomials.
\newblock {\em Proc. Steklov Inst. Math.}, 280:S65--S70, 2013.

\bibitem{Poonen}
Bjorn Poonen.
\newblock The classification of rational preperiodic points of quadratic
  polynomials over: a refined conjecture.
\newblock {\em Mathematische Zeitschrift}, 228(1):11--29, 1998.

\bibitem{schoof2010catalan}
Ren{\'e} Schoof.
\newblock {\em Catalan's conjecture}.
\newblock Springer Science \& Business Media, 2010.

\bibitem{SilvermanBook}
Joseph~H. Silverman.
\newblock {\em The arithmetic of dynamical systems}, volume 241 of {\em
  Graduate Texts in Mathematics}.
\newblock Springer, New York, 2007.

\bibitem{Stoll}
Michael Stoll.
\newblock Galois groups over {${\bf Q}$} of some iterated polynomials.
\newblock {\em Arch. Math. (Basel)}, 59(3):239--244, 1992.

\bibitem{stoll2008rational}
Michael Stoll.
\newblock Rational 6-cycles under iteration of quadratic polynomials.
\newblock {\em LMS Journal of Computation and Mathematics}, 11:367--380, 2008.

\bibitem{stoll2011rational}
Michael Stoll.
\newblock Rational points on curves.
\newblock {\em Journal de th{\'e}orie des nombres de Bordeaux}, 23(1):257--277,
  2011.

\bibitem{FLT}
Andrew Wiles.
\newblock Modular elliptic curves and fermat's last theorem.
\newblock {\em Annals of mathematics}, 141(3):443--551, 1995.

\end{thebibliography}
\end{document}